\documentclass[11pt]{amsart}
\usepackage{mathrsfs}
\usepackage{amssymb}
\usepackage[all,pdf]{xy}
\usepackage{titletoc}
\usepackage{amscd}
\usepackage{amsmath, amssymb}
\usepackage{amsfonts}
\usepackage[colorlinks,linkcolor=blue,citecolor=blue, pdfstartview=FitH]{hyperref}

\numberwithin{equation}{section}

\theoremstyle{plain}
\newtheorem{thm}{Theorem}[section]
\newtheorem{theorem}[thm]{Theorem}
\newtheorem{lemma}[thm]{Lemma}
\newtheorem{corollary}[thm]{Corollary}
\newtheorem{proposition}[thm]{Proposition}

\theoremstyle{definition}

\newtheorem{remark}[thm]{Remark}

\newtheorem{definition}[thm]{Definition}

\newtheorem{example}[thm]{Example}

\newtheorem{defn-thm}[thm]{Definition-Theorem}

\newcommand{\sC}{{\mathcal C}}

\newcommand{\sT}{{\mathcal T}}

\newcommand{\C}{{\mathbb C}}

\newcommand{\F}{{\mathbb F}}

\renewcommand{\P}{{\mathbb P}}

\newcommand{\T}{{\mathbb T}}

\newcommand{\Z}{{\mathbb Z}}

\newcommand{\Img}{{ \mathrm{Im}}}
\newcommand{\Hom}{{\mathrm{Hom}}}
\newcommand{\Ker}{{\mathrm{Ker}}}
\newcommand{\Coker}{{\mathrm{Coker}}}

\newcommand{\Ext}{\mathrm{Ext}}
\newcommand{\Gr}{{\mathrm{Gr}}}

\newcommand{\btheorem}{\begin{theorem}}
\newcommand{\etheorem}{\end{theorem}}
\newcommand{\bproposition}{\begin{proposition}}
\newcommand{\eproposition}{\end{proposition}}
\newcommand{\bdefinition}{\begin{definition}}
\newcommand{\edefinition}{\end{definition}}
\newcommand{\bcorollary}{\begin{corollary}}
\newcommand{\ecorollary}{\end{corollary}}
\newcommand{\bproof}{\begin{proof}}
\newcommand{\eproof}{\end{proof}}
\newcommand{\bremark}{\begin{remark}}
\newcommand{\eremark}{\end{remark}}
\newcommand{\eexample}{\end{example}}
\newcommand{\bexample}{\begin{example}}

\newcommand{\elemma}{\end{lemma}}
\newcommand{\blemma}{\begin{lemma}}

\newcommand{\lrw}{\longrightarrow}

\usepackage{tikz-cd}
\usepackage[all]{xy}
\usepackage{fancyhdr}
\begin{document}

\title{Motivic cluster multiplication formulas in 2-Calabi-Yau categories}
\author{Jie Xiao, Fan Xu, Fang Yang*}
\address{School of Mathematical Sciences\\
Beijing Normal University\\
Beijing 100875, P.~R.~China}
\email{jxiao@bnu.edu.cn(J. Xiao)}
\address{Department of Mathematical Sciences\\
Tsinghua University\\
Beijing 100084, P.~R.~China}
\email{fanxu@mail.tsinghua.edu.cn(F. Xu)}
\address{Department of Mathematical Sciences\\
Tsinghua University\\
Beijing 100084, P.~R.~China}
\email{yangfang19@mails.tsinghua.edu.cn(F. Yang)}
\thanks{$*$~Corresponding author.}

\begin{abstract}
  We introduce a notion of motivic cluster characters via virtual Poincar\'{e} polynomials, and prove a motivic version of multiplication formulas obtained by Chen-Xiao-Xu for weighted quantum cluster characters associated to a 2-Calabi-Yau triangulated category $\sC$ with a cluster tilting object. Furthermore, a refined form of this formula is also given. When $\sC$ is the cluster category of an acyclic quiver, our certain refined multiplication formula is a motivic version of the multiplication formula in \cite{Chen2023}.
\end{abstract}

\maketitle
\setcounter{tocdepth}{1} \tableofcontents

\section{Introduction}\label{sec0}
  It is known that cluster algebras are closely related to representation theory of quivers via cluster categories (over $\C$) and cluster characters. Cluster category admits cluster tilting objects and has the $2$-Calabi-Yau property. Cluster characters were firstly introduced by Caldero and Chapoton for Dynkin quivers \cite{Caldero2006}, then for any acyclic quiver by Caldero and Keller \cite{Caldero2006a}, and finally generalized by Palu to any 2-Calabi-Yau triangulated category $\sC$ with a cluster tilting object $T$ \cite{Palu2008}. For $M\in \sC$, the cluster character of $M$ is given as follows ,
      $$CC_T(M)=\sum_{e} \chi(\Gr_{e}(\Hom_{\sC}(T,M)))X^{p(M,e)}.$$
  where $\chi$ is the Euler-Poincar\'{e} characteristic (see Section \ref{sec2.2} for more details).
  
  Cluster characters build a bridge between cluster algebras and cluster categories. The essential point lies in mutations of quivers with potentials and cluster multiplication formulas. For any $M$, $N\in \sC$ such that $\dim \Ext^1_{\sC}(M,N)=1$, we have 
    \begin{equation}\label{in1}
    CC_T(M)\cdot CC_T(N)=CC_T(L)+ CC_T(L'),
    \end{equation}
  where $L$ and $L'$ are the middle terms of exchange triangles of $M$ and $N$. The formula (\ref{in1}) was firstly proved by by Caldero-Chapoton for Dynkin quivers \cite{Caldero2006}. Since then, cluster multiplication formulas have been studied extensively, e.g. Caldero-Keller \cite{Caldero2008}, Xiao-Xu \cite{Xiao2009, Xiao2010a}, Xu \cite{Xu2010}, and Palu to any 2-Calabi-Yau triangulated category with a cluster tilting object \cite{Palu2012}. For any $M$, $N\in \sC$, Palu showed that
  \begin{equation}\label{in2}
    \begin{aligned}
    &\chi(\mathbb{P}\Ext^1_{\sC}(M,N))CC_T(M)\cdot CC_T(N)\\
    =&\sum_{[L]}(\chi(\mathbb{P}\Ext^1_{\sC}(M,N)_L)+\chi(\mathbb{P}\Ext^1_{\sC}(N,M)_L))CC_T(L).
    \end{aligned}
  \end{equation}
  The proof heavily depends on the 2-Calabi-Yau property. A refined version of this multiplication formula recently has been shown by Keller-Plamondon-Qin in \cite{Keller2023}. For a non-zero vector subspace of $\Ext^1_{\sC}(M,N)$, they proved the following identity (see Section \ref{secref} for more details):
  \begin{equation}\label{in3}
    \begin{aligned}
    &\chi(\P V)CC_T(M)CC_T(N)\\
    =&\sum_{E\not\cong M\oplus N} (\chi(\P V_E)+ \chi(\P(\Ext^1_{\sC}(N,M)\setminus \P(V^{\perp}))_E)CC_T(E).
    \end{aligned}
  \end{equation}
  
  Quantum cluster algebra is a quantization of cluster algebra introduced by Berenstein and Zelevinsky \cite{Berenstein2005}. In the representation theory side, quantum cluster characters are quantum analogue of cluster characters, constructed by Rupel \cite{Rupel2011} for acyclic quivers. Let $Q$ be an acyclic quiver and $k=\F_q$ a finite field. For $M\in \mathrm{mod}(kQ)$, 
      $$\tilde{X}_{M}:=\int_{e} q^{-\frac{1}{2}\langle e,m-e\rangle} |\mathrm{Gr}_e(M)| X^{-e^*-^*(m-e)}.$$
  Qin used an alternative definition for quantum cluster characters to show the quantized version of multiplication formula (\ref{in1}) for $\tilde{X}_M$ and $\tilde{X}_N$. This quantized multiplication formula was then generalized by Ding-Xu for affine cluster algebras \cite{Ding2012a}. The higher-dimensional multiplication formula of quantum cluster characters for any acyclic quiver has been shown by Chen-Ding-Zhang in \cite{Chen2023}.
  \begin{equation*}
    \begin{aligned}
     &(q^{[M,N]^1}-1) \tilde{X}_M\cdot \tilde{X}_N= q^{\frac{1}{2}\Lambda(m^*,n^*)}\sum_{E\not\cong M\oplus N} |\Ext^1_{kQ}(M,N)_E|\tilde{X}_E\\
     &+ \sum_{D\not\cong N,A,I}q^{\frac{1}{2}\Lambda((m-a)^*,(n+a)^*)+\frac{1}{2}<m-a,n>}|_D\Hom_{kQ}(N,\tau M)_{\tau A\oplus I}| \tilde{X}_A\cdot\tilde{X}_{D\oplus I[-1]}.
    \end{aligned}
  \end{equation*}
  Xu-Yang in \cite{Xu2022} gave an explicit formulas for quantum cluster characters from cluster categories of weighted projective lines and established a similar formula as above in the categories of coherent sheaves on weighted projective lines. To the authors' knowledge, the explicit Laurent expansion formulas for quantum cluster characters associated to  any 2-Calabi-Yau triangulated category with a cluster tilting object have still been unknown. Huang obtained a Laurent expansion formula for the quantum cluster algebras from unpunctured triangulated surfaces \cite{Huang2022} and unpunctured orbifolds with arbitrary coefficients and quantization \cite{Huang2022a}. However, the quantum version of the multiplication formula (\ref{in2}) has been obtained by Chen-Xiao-Xu by utilizing weighted quantum cluster characters \cite{Chen2021a}, though their weighted cluster characters depend on the finite ground field. 
    
  A question is whether we could find a generic version of multiplication formulas generalizing all (quantized or not) multiplication formulas mentioned above for any 2-Calabi-Yau triangulated category $\sC$ with a cluster tilting object $T$. The present paper replies this question affirmatively. Note that all multiplication formulas mentioned above heavily depend on 2-Calabi-Yau property. Indeed, For any $M$, $N\in \sC$, set $F:=\Hom_{\sC}(T,-)$ and
     $$Z_1:=\{(\epsilon,L_0)|\ \epsilon\in\Ext^1_{\sC}(M,N),L_0\in \bigsqcup_{g}\Gr_g(F\mathrm{mt}\epsilon)\},$$
  and $\rm{Pr}_1: Z_1\to \Ext^1_{\sC}(M,N)$. Set 
     $$Z_2:=\{(\eta,L'_0)|\ \eta\in\Ext^1_{\sC}(N,M),L'_0\in \bigsqcup_{g}\Gr_g(F\mathrm{mt}\eta)\},$$
  and $\rm{Pr}'_1:Z_2\to \Ext^1_{\sC}(N,M)$.  We have a map (see Section \ref{motsec})
  \begin{equation*}
    \begin{aligned}
    &\Psi:=\Psi_1\sqcup\Psi_2: Z_1\sqcup Z_2\lrw \bigsqcup_{e+f=g}\Gr_e(FM)\times \Gr_f(FN),\\
    &\Psi_1: (\epsilon,L_0)\mapsto \psi^{\epsilon}(L_0),\qquad \Psi_2: (\eta,L'_0)\mapsto \psi^{\eta}(L'_0).
    \end{aligned}
  \end{equation*}
  The orthogonal complement of $\rm{Pr}_1(\Psi_1^{-1}(M_0,N_0))$ is precisely $\rm{Pr}'_1(\Psi_2^{-1}(M_0,N_0))$ for any $(M_0,N_0)$. Moreover the fiber of $\Psi$ is an affine space, which means that we can explicitly calculate the size of fibers. For the refined case, we focus on a non-zero vector subspace $V$ instead of the whole space $\Ext_{\sC}^1(M,N)$. In this case, we need take a subspace or subset $W$ of $\Ext^1_{\sC}(N,M)$ such that 
  $$\dim W\cap \rm{Pr}'_1(\Psi_2^{-1}(M_0,N_0))+\dim V\cap \rm{Pr}_1(\Psi_1^{-1}(M_0,N_0))$$
  is fixed for any $(M_0,N_0)$. In \cite{Keller2023}, $W$ is chosen to be the subset $\Ext^1_{\sC}(N,M)\setminus V^{\perp}$ . 

  To construct such a generic version of multiplication formulas, we use motivic weighted cluster characters, which become cluster characters after evaluating $t$  to $-1$. More precisely, for $M\in \sC$, we define the motivic cluster character of $M$ to be 
      $$X^T_M=\sum_{e}\Upsilon(\mathrm{Gr}_e(\Hom_{\sC}(T,M))) X^{p(M,e)},$$
  where $\Upsilon(Z)$ is the virtual Poincar\'{e} polynomial of locally closed set $Z$. Then for a $\Z$-valued weight function $f$ and $\epsilon\in \Ext^1_{\sC}(M_1,M_2)_M$, we could similarly define a motivic weighted cluster character $f(\epsilon,-)*X^T_M$.

  The paper is organized as follows. In Section \ref{motsec} we construct a motivic version of multiplication formulas for any two weighted motivic cluster characters. Then Section \ref{secref} is devoted to establish a refined motivic version of multiplication formulas, which generalizes the formula given in Section \ref{motsec}, also recovers the Keller-Plamondon-Qin's result. Finally in Section \ref{SecHer}, we establish several specific refined motivic multiplication formulas in hereditary case, which is a motivic version of Chen-Ding-Zhang formulas in \cite{Chen2023}.

\section{2-Calabi-Yau categories and motivic  multiplication formulas}\label{motsec}
  \subsection{Exact structure}\label{sec2.1}
    Let $\sC$ be a $\Hom$-finite, Krull-Schimdt triangulated category over $\C$ with a suspension functor $\Sigma$. We say that $\sC$ is  $2$-$Calabi$-$Yau$ if for any objects $M$ and $N$ of $\sC$, it is equipped with a bilinear form
        $$\beta_{MN}:\Hom_{\sC}(M,\Sigma N)\otimes \Hom_{\sC}(N,\Sigma M)\lrw \C,$$
    which is non-degenerate and bifunctorial.
    Let $T$ be a cluster-tilting object in $\sC$ and $B$ be the endomorphism algebra $\mathrm{End}_{\sC}(T)$ of $T$, then there is a functor
      $$F:=\Hom_{\sC}(T,-):\sC\lrw \mathrm{mod}B,\ X \mapsto \Hom_{\sC}(T,X),$$
    which induces an equivalence of categories $\sC/(\Sigma T)\stackrel{\cong}\lrw \mathrm{mod}B$. 
    
    Following from Lemma \cite[Lemma 1.3]{Palu2008}, the skew-symmetric form $\langle-,-\rangle_a$ of the partial Euler form $\langle-,-\rangle$ is a bilinear form on the Grothendieck group $K_0(\mathrm{mod}B)$, where $\langle-,-\rangle$ is given by 
       $$\langle [M],[N]\rangle=\dim \Hom_B(M,N)-\dim \mathrm{Ext}_B(M,N),$$
    for $[M]$, $[N] \in K_0(\mathrm{mod}B)$. For $X\in \sC$, there are two triangles \cite[Section 2.1]{Palu2008},
       $$T_2\to T_1\to X\to \Sigma T_2\qquad \text{and}\qquad \Sigma T_4\to X\to \Sigma^2 T_3\to \Sigma^2 T_4.$$
    with all $T_i$ in $\mathrm{add}T$. Then we can define 
    index and coindex of $X$ by
        $$\rm{Ind}X=[FT_1]-[FT_2]\qquad \text{and} \qquad \rm{coind}X=[FT_3]-[FT_4].$$
    
    A given $\epsilon\in \Hom_{\sC}(M,\Sigma N)$ induces a triangle 
       $$N\stackrel{i}\lrw L\stackrel{p}\lrw M\stackrel{\epsilon}\lrw \Sigma N$$
    in $\sC$. We call $L$ the middle term of $\epsilon$ and denote by $\mathrm{mt}\epsilon$. Applying the functor $F$, we get an exact sequence in $\mathrm{mod}B$
       $$FN\stackrel{Fi}\lrw FL\stackrel{Fp}\lrw FM\stackrel{F\epsilon}\lrw F(\Sigma N).$$
    For $g\in K_0(\mathrm{mod}B)$, denote by $\Gr_g(FL)$ the set of all submodules of $FL$ with dimension vector $g$. Consider the following map
      $$\psi^{\epsilon}_g: \Gr_{g}(FL)\lrw \bigsqcup_{e+f=g+[\Ker(Fi)]} \Gr_e(FM)\times Gr_{f}(FN),\ L_0\mapsto (Fp(L_0),Fi^{-1}(L_0)).$$  
    Note that $\psi^{\epsilon}_g$ may not be surjective in general, we should figure out the image $\Img(\psi)$ and fibers of $\psi^{\epsilon}_g$
      $$ \mathrm{Gr}_{M_0,N_0}^{\epsilon}(FL):=\{L_0\subseteq FL|\ Fi^{-1}(L_0)=N_0,\ Fp(L_0)=M_0\}.$$
    Similarly for $\eta:M\stackrel{i'}\lrw L'\stackrel{p'}\lrw N\stackrel{\eta}\lrw \Sigma M\in \Ext^1_{\sC}(N,M)$, define
      $$\psi^{\eta}_g: \Gr_{g}(FL')\lrw \bigsqcup_{\substack{e+f=\\g+[\Ker(Fi')]}} \Gr_f(FN)\times Gr_{e}(FM),\ L_0\mapsto (Fp'(L_0),(Fi')^{-1}(L_0)),$$  
    and fibers
      $$ \mathrm{Gr}_{N_0,M_0}^{\eta}(FL'):=\{L_0\subseteq FL'|\ (Fi')^{-1}(L_0)=M_0,\ Fp(L'_0)=N_0\}.$$
    Note that $\rm{Gr}_g(FL)$ (resp.  $\rm{Gr}_g(FL')$) admits a stratification of locally closed subsets such that $\psi_g^{\epsilon}$ (resp. $\psi^{\eta}_{g}$) restricted to each strata is a morphism of varieties. Thus the pull-back of constructible subsets along $\psi_g^{\epsilon}$ (resp. $\psi_g^{\eta}$) are still constructible. In the sequel, write $\psi^{\epsilon}:=\sqcup_{g}\psi_g^{\epsilon}$ and  $\psi^{\eta}:=\sqcup_{g}\psi_g^{\eta}$. 

    Let $\widetilde{M_0}\stackrel{\iota_m}\to M$ and $\widetilde{N_0}\stackrel{\iota_N}\to N$ be lifts of $M_0\subseteq M$, $N_0\subseteq N$ along $F$. Following from \cite[Section 4]{Palu2008} and \cite[Section 4.2]{Chen2021a}, we have two linear maps:
    \begin{equation*}
      \begin{aligned}
      \alpha_{M_0,N_0}\colon&\Hom_{\mathcal{C}}(\Sigma^{-1}M,\tilde{N}_0)\oplus\Hom_{\mathcal{C}}(\Sigma^{-1}M,N)\\
      &\longrightarrow\Hom_{\mathcal{C}/(T)}(\Sigma^{-1}\tilde{M}_0,\tilde{N}_0)\oplus\Hom_{\mathcal{C}}(\Sigma^{-1}\tilde{M}_0,N)\oplus\Hom_{\mathcal{C}/(\Sigma T)}(\Sigma^{-1}M,N)\\
      &(a,b) \longmapsto (a\circ\Sigma^{-1}\iota_M,\iota_N\circ a\circ\Sigma^{-1}\iota_M-b\circ\Sigma^{-1}\iota_M,\iota_N\circ a-b)
      \end{aligned}
      \end{equation*}
      and
    \begin{equation*}
      \begin{aligned}
      \alpha'_{N_0,M_0}\colon&\Hom_{\Sigma T}(\tilde{N}_0,\Sigma\tilde{M}_0)\oplus\Hom_{\mathcal{C}}(N,\Sigma \tilde{M}_0)\oplus\Hom_{\Sigma^2T}(N,\Sigma M)\\
      &\longrightarrow\Hom_{\mathcal{C}}(\tilde{N}_0,\Sigma M)\oplus\Hom_{\mathcal{C}}(N,\Sigma M)\\
      &(a,b,c)\longmapsto (\Sigma\iota_M\circ a+c\circ\iota_N+\Sigma\iota_M\circ b\circ\iota_N,-c-\Sigma\iota_M\circ b).
      \end{aligned}
    \end{equation*}
    
    \begin{lemma}[{\cite[Section 4.2]{Chen2021a}}]\label{fiberlem}
    Fix $(M_0,N_0)\in \Gr_e(FM)\times \Gr_f(FN)$.\\
    (i) For $\epsilon\in \mathrm{Ext}_{\sC}^1(M,N)$, $(M_0,N_0)\in \Img\psi^{\epsilon}$ if and only if 
       $$\epsilon\in \Sigma p\Ker(\alpha_{M_0,N_0}),$$
    where $p$ is the natural projection to $Ext^1_{\sC}(M,N)$. Furthermore in this case, the fiber of $(M_0,N_0)$ is an affine space and satisfies
       $$\mathrm{Gr}_{M_0,N_0}^{\epsilon}(F(\mathrm{mt}\epsilon))\cong \Hom_{B}(M_0,FN/N_0).$$ 
    (ii) For $\eta\in \mathrm{Ext}^1_{\sC}(N,M)$, $(N_0,M_0)\in \Img\psi^{\eta}$ if and only if 
       $$\eta\in \mathrm{Im}(\alpha'_{N_0,M_0})\cap \mathrm{Ext}^1_{\sC}(N,M).$$
    Furthermore in this case, the fiber of $(N_0,M_0)$  is an affine space and satisfies
        $$\mathrm{Gr}_{N_0,M_0}^{\eta}(F(\mathrm{mt}\eta))\cong \Hom_{B}(N_0,FM/M_0).$$
    (iii) $\Sigma\Ker(\alpha_{M_0,N_0})=(\mathrm{Im}\alpha'_{N_0,M_0})^{\perp}$. Moreover, 
      $$\dim \Sigma p\Ker\alpha_{M_0,N_0} +\dim (\mathrm{Im}\alpha'_{N_0,M_0}\cap \Ext^1_{\sC}(N,M))=\dim \Ext^1_{\sC}(M,N).$$ 
    \end{lemma}
  It can be deduced that $\psi_g^{\epsilon}$ and $\psi_g^{\eta}$ are locally of finite representation by the lemma above. Thus  $\psi_g^{\epsilon}$ and $\psi_g^{\eta}$ preserve constructible subsets.
  \begin{lemma}
    For $\epsilon\in \Ext^1_{\sC}(M,N)$ and $\eta\in \Ext^1_{\sC}(N,M)$, subsets $\Img\psi_g^{\epsilon}$ and $\Img\psi_g^{\eta}$ are constructible.
  \end{lemma}

  Fix $\epsilon\in \Ext^1_{\sC}(M,N)$ as above, define a $\Z$-valued function as follows:
     $$\begin{tikzcd}
      \mathrm{Gr}_g(FL)\arrow[r,"\psi_g^{\epsilon}"]\arrow[rr,dotted,bend left=20,"\kappa_g^{\epsilon}"]
     &\bigsqcup\limits_{\substack{e+f=\\g+[\Ker(Fi)]}} \mathrm{Gr}_e(FM)\times \mathrm{Gr}_f(FN)\arrow[r,"d_{M,N}"]  &\Z,
     \end{tikzcd}$$
  where $d_{M,N}$ sends $(M_0,N_0)$ to $\dim_{\C}\Hom_B(M_0,FN/N_0)$, which is an upper semicontinuous function following from the lemma below.

  \begin{lemma}[{\cite[Lemma 4.3]{Crawley2002}}]
    Let $\mathrm{mod}_B^e$ be the subcategory of $\mathrm{mod}_B$ of modules with dimension vector $e$. The function 
      $$\mathrm{mod}_B^{e}\times \mathrm{mod}_B^f\lrw \Z$$
    sending a pair $(M_1,M_2)$ to the dimension of $\Hom_B(M_1,M_2)$ is upper semicontinuous.
  \end{lemma}

  As a consequence, $d_{M,N}^{-1}(c)$ and $(\kappa_g^{\epsilon})^{-1}(c)$ are locally closed subsets,  and $\Img\psi_g^\epsilon$ has a partition
     $$\Img\psi_g^{\epsilon}=\bigsqcup_{c\in \Z}d_{M,N}^{-1}(c)\cap \Img\psi_g^{\epsilon}=:\bigsqcup_{c\in \Z} Z_{g,c}^{\epsilon}.$$
  Also we have
     $$\mathrm{Gr}_g(FL)\cong \bigsqcup_{c\in \Z} (\kappa_g^{\epsilon})^{-1}(c)=:\bigsqcup_{c\in \Z} F_{g,c}^{\epsilon},$$
  where $F_{g,c}^{\epsilon}$ is an affine fibration over $Z_{g,c}^{\epsilon}$ of rank $c$ by Lemma \ref{fiberlem}(i), .
  Similarly for $\eta\in \Ext^1_{\sC}(N,M)$, we have
     $$\mathrm{Gr}_g(F\mathrm{mt}(\eta))\cong \bigsqcup_{c\in \Z} F_{g,c}^{\eta},$$
  where $F_{g,c}^{\eta}$ is an affine fibration over $Z_{g,c}^{\eta}$ of rank $c$ and $Z_{g,c}^{\eta}=d_{N,M}^{-1}(c)\cap \Img\psi_{g}^{\eta}$.

\subsection{Virtual Poincar\'{e} polynomials}
  For any locally closed subset $Z$ of a smooth projective variety over $\C$, we can assign a  virtual Poincar\'{e} polynomial $\Upsilon(Z)$ to $Z$, such that:
  \begin{itemize}
    \item [(i)]  For any smooth projective variety X, $\Upsilon(X)=\sum_{i=1}^{2\dim X} \dim \text{H}^i_c(X,\C)t^i$;
    \item [(ii)] For $X=X_1\sqcup X_2$, then $\Upsilon(X)=\Upsilon(X_1)+\Upsilon(X_2)$;
    \item [(iii)] For any affine fibration $f:X\to Y$ with fiber $Z$, then $\Upsilon(X)=\Upsilon(Y)\Upsilon(Z)$.
  \end{itemize}
  This assignment is unique (see \cite[Theorem 2.14]{Joyce2007a}) and \cite[Section 5.3]{Bridgeland2017}). By the first condition, we can deduce that $\Upsilon(Z)|_{t=-1}$ is precisely the Euler-Poincar\'{e} characteristic $\chi(Z)$. According the third condition above, we have
  
  \begin{lemma}\label{lem2.5}
    For $\epsilon\in \Ext^1_{\sC}(M,N)$ and $\eta\in \Ext^1_{\sC}(N,M)$, there are two identities of virtual Poincar\'{e} polynomials 
        $$\Upsilon(\mathrm{Gr}_g(F(\mathrm{mt}\epsilon)))=\sum_{c\in \Z}t^{2c}\Upsilon(Z_{g,c}^{\epsilon})$$
    and  
        $$\Upsilon(\mathrm{Gr}_g(F(\mathrm{mt}\eta)))=\sum_{c\in \Z}t^{2c}\Upsilon(Z_{g,c}^{\eta})$$
       
  \end{lemma}

\subsection{Motivic weighted cluster characters}\label{sec2.2}
  Write $K_0(\mathrm{mod}B)\cong \Z^n$. Let $\Lambda(-,-): K_0(\mathrm{mod}B)\times K_0(\mathrm{mod}B)\to \Z$ be a skew-symmetric form. The quantum torus $\sT_{\Lambda}$ is defined to be 
      $$\sT_{\Lambda}:=\C(t)[X^{\alpha}|\ \alpha\in \Z^n],$$
  with a twisted multiplication given by
      $$X^e \cdot X^f=t^{\Lambda(e,f)}X^{e+f}.$$
  
  \begin{definition}
  For each object $L\in \sC$, the motivic cluster character of $L$ is defined to be
       $$X_L:=\int_g \Upsilon(\mathrm{Gr}_g(FL)) X^{p(L,g)}=\int_g\int_{(M_0,N_0)\in \mathrm{Gr}_g(FL)} X^{p(L,g)},$$
  where $p(L,g)=-\mathrm{coind}_T L+C\cdot g$ and $C$ is the matrix whose $(i,j)$-term is $\langle[S_i],[S_j]\rangle_a$. Here $S_i$ are simple $B$-modules.
  \end{definition}

  Recall that we have partition on $\mathrm{Gr}_g(FL)$ associated to $\epsilon$ in Section \ref{sec2.1}:
      $$\mathrm{Gr}_g(FL)=\bigsqcup_{c\in \Z} F_{g,c}^{\epsilon}$$
  where $F_{g,c}^{\epsilon}$ is an affine fibration over $Z_{g,c}^{\epsilon}$ of rank c.
  Using Lemma \ref{lem2.5}, it can be seen that for $\epsilon\in \Ext^1_{\sC}(M,N)$ with middle term $L$, we have
      $$X_L=\int_g\int_{c\in \Z} t^{2c}\Upsilon(Z^{\epsilon}_{g,c})X^{p(L,g)}=\int_g\int_{c}\int_{(M_0,N_0)\in Z^{\epsilon}_{g,c}} t^{2c} X^{p(L,g)}.$$

  For any $M$, $N\in \sC$, set $MG(M,N)$ to be the following constructible set
      $$ \{(\epsilon,M_0,N_0)\in \Ext^1_{\sC}(M,N)\times \rm{Gr}(FM)\times \rm{Gr}(FN)|\ (M_0,N_0)\in \Img\psi^{\epsilon}\}.$$
  \begin{definition}
  A weight function associated to $(M,N)$ is a $\Z$-valued constructible function defined on $MG(M,N)$. 
  \end{definition}
  Thus sending $(\epsilon,M_0,N_0)\in \bigsqcup_g Z_{g,c}^{\epsilon}$ to $-2c$ gives rise to a weight function $f_{M,N}$. Similarly, we have another weight function $f_{N,M}$, which sends $(\eta,N_0,M_0)\in \bigsqcup_g Z_{g,c}^{\eta}$ to $-2c$.
      
  For $\epsilon\in \Ext^1_{\sC}(M,N)$, denote by $f(\epsilon,-)$ be the weight function $f$ restricted to $MG(M,N)\cap (\{\epsilon\}\times \rm{Gr}(FM)\times \rm{Gr}(FN))$.
  \begin{definition}
  The motivic weighted cluster character $f(\epsilon,-)*X_L$ of $L$ by $f(\epsilon,-)$ with $\rm{mt}\epsilon=L$ is given by
  \begin{equation*}
    \begin{aligned}
    f(\epsilon,-)*X_L:=&\int_{g}\int_{L_0\in \Gr_{g}(FL)}t^{f(\epsilon,\psi^{\epsilon}_{g}(L_0))}X^{p(L,g)}\\
    =&\int_g\int_c\int_{(M_0,N_0)\in Z^{\epsilon}_{g,c}} t^{2c}t^{f(\epsilon,M_0,N_0)}X^{p(L,g)}.
    \end{aligned}
  \end{equation*}
  \end{definition}
  For $\lambda\in \C^{\times}$, note that $\Img\psi_g^{\epsilon}=\Img\psi_{g}^{\lambda\epsilon}$ by Lemma \ref{fiberlem} and $\mathrm{mt}\epsilon\cong \mathrm{mt}(\lambda\epsilon)$. For the projectivization of $\Ext^1_{\sC}(M,N)$, we assume that
      $$f(\epsilon,-)=f(\lambda\epsilon,-),$$
  for all $\lambda\in \C^{\times}$. With the assumption, we can see that
      $$f(\epsilon,-)*X_L=f(\lambda\epsilon,-)*X_L.$$
  For $\epsilon\in \Ext^1_{\sC}(M,N)$ such that $L=\mathrm{mt}\epsilon$, if we take  $f_{M,N}$ to be the weight function, then the motivic cluster character of $L$ weighted by $f_{M,N}$ is
    \begin{equation*}
    \begin{aligned}
      f_{M,N}*(\epsilon,-)X_{\rm{mt}\epsilon}
      =&\int_g\int_{L_0\in \Gr_g(FL)} t^{f(\epsilon,\psi_g^{\epsilon}(L_0))}X^{p(\rm{mt}\epsilon,g)}\\
      =&\int_g\int_c\int_{(M_0,N_0)\in Z^{\epsilon}_{g,c}}
      t^{2c}t^{f(\epsilon,M_0,N_0)}X^{p(\rm{mt}\epsilon,g)} \\
      =&\int_g\int_c \Upsilon(Z^{\epsilon}_{g,c}) X^{p(\rm{mt}\epsilon,g)}.  
    \end{aligned}
    \end{equation*}

\subsection{Motivic weighted multiplication formulas}
  In this subsection, we will show  multiplication formulas for motivic weighted cluster characters. Fix $M$, $N\in \sC$ such that $\Ext^1_{\sC}(M,N)\neq 0$. Set $g^{+}$ to be the weight function sending $(\epsilon,M_0,N_0)\in MG(M,N)$ to $d_{M_0,N_0}$=$\dim \Img\alpha'_{N_0,M_0}\cap \Ext^1_{\sC}(N,M)$, and $g^-$ to be the weight function sending $(\eta,N_0,M_0)\in MG(N,M)$ to $d'_{N_0,M_0}$= $\dim \Sigma p\Ker\alpha_{M_0,N_0}$. Note that pairs $(g^+(\epsilon,-),0)$ and $(0,g^-(\eta,-))$ are the pointwise balanced pairs defined in \cite[Section 4]{Chen2021a}. Recall that $(h^+,h^-)$ is called a $pointwise$ $balanced$ $pair$ if
     \begin{align*}
      \Upsilon(\P\Ext^1_{\sC}(M,N))=&t^{h^+(M_0,N_0)}\Upsilon(\P(\Sigma p\Ker\alpha_{M_0,N_0}))+\\
      &\qquad t^{h^-(N_0,M_0)}\Upsilon(\P(\Img\alpha'_{N_0,M_0}\cap \Ext^1_{\sC}(N,M)))
     \end{align*}
  for any $(M_0,N_0)\in \Gr(FM)\times \Gr(FN)$.
  Recall that we have defined two weight functions $f_{M,N}$ and $f_{N,M}$ putting weights $-2c$ on $\{\epsilon\}\times \bigsqcup_g Z^{\epsilon}_{g,c}$ for $\epsilon\in \Ext^1_{\sC}(M,N)$ and $\{\eta\}\times\bigsqcup_g Z^{\eta}_{g,c}$ for $\eta\in \Ext^1_{\sC}(N,M)$  respectively. Define another two weight functions associated to $(M,N)$ and $(N,M)$ by 
     $$\bar{f}_{M,N}(\epsilon,M_0,N_0):=f_{M,N}(\epsilon,M_0,N_0)+ \Lambda(p(M,[M_0]),p(N,[N_0])),$$
  and  
     $$\bar{f}_{N,M}(\eta,N_0,M_0):=f_{M,N}(\eta,N_0,M_0)+ \Lambda(p(M,[M_0]),p(N,[N_0])).$$

  Let $g(\xi,-)*X_M$ and $g'(\xi',-)*X_N$ be two motivic weighted cluster characters. Define two weight functions of $g(\xi,-)$ and $g'(\xi',-)$ as follows:
    \begin{align*}
       &\T^{\xi,\xi'}_{g,g'}: MG(M,N)\lrw \Z\\
       &\qquad\ \ (\epsilon,M_0,N_0)\mapsto g(\xi,\psi^{\xi}(M_0))+g'(\xi',\psi^{\xi'}(N_0)),
    \end{align*}
  and 
    \begin{align*}
      &\hat{\T}^{\xi,\xi'}_{g,g'}: MG(N,M)\lrw \Z\\
      &\qquad\ \ (\eta,N_0,M_0)\mapsto g(\xi,\psi^{\xi}(M_0))+g'(\xi',\psi^{\xi'}(N_0)).
    \end{align*}

  \begin{lemma}[{\cite[Lemma 5.1]{Palu2008}}]\label{lemp}
  For $\epsilon\in \Ext^1_{\sC}(M,N)_L$, if $\psi_g^{\epsilon}(L_0)=(M_0,N_0)\in \Gr_e(FM)\times \Gr_f(FN)$, then
     $$p(L,g)=p(M,e)+p(N,f).$$
  \end{lemma}
 
  Note that for the split extension $0_{M\oplus N}\in \Ext^1_{\sC}(M, N)_{M\oplus N}$, the map
     $$\psi_g^0: \bigsqcup_{g}\Gr_g(F(M\oplus N))\lrw \bigsqcup_{e+f=g} \Gr_e(FM)\times \Gr_f(FN)$$
  is an affine fibration with fiber isomorphic to $\Hom_{\sC}(M_0,FN/N_0)$. Now we can state the following 
  \begin{proposition}\label{prop2.11}
  For any $M$, $N\in \sC$, let $g(\xi,-)*X_M$ and $g'(\xi',-)*X_N$ be two motivic weighted cluster characters. Then
  \begin{equation}
    \begin{aligned}
      g(\xi,-)*X_M\cdot g'(\xi',-)*X_N=\bar{f}_{M,N}(0,-)*\T_{g,g'}^{\xi,\xi'}(0,-)* X_{M\oplus N}
    \end{aligned}
  \end{equation}

  \begin{proof}
  To simplify notations, denote by $\T:=\T_{g,g'}^{\xi,\xi'}$, we have that
   \begin{equation*}
    \begin{aligned}
      \mathrm{RHS}=&\int_g\int_c\int_{(M_0,N_0)\in Z^{0}_{g,c}} t^{2c}\cdot t^{f(0,M_0,N_0)+\Lambda(p(M,[M_0]),p(N,[N_0]))}\cdot t^{\T(0,M_0,N_0)}\cdot X^{p(M\oplus N,g)}\\
      =&\int_g\int_c\int_{(M_0,N_0)\in Z^0_{g,c}} t^{\T(0,M_0,N_0)}\cdot X^{p(M,[M_0])}\cdot X^{p(N,[N_0])}\\
      =&\int_{e,f}\int_{(M_0,N_0)\in \Gr_e(FM)\times \Gr_f(FN)}  t^{g(\xi,\psi^{\xi}_e(M_0))}\cdot t^{g'(\xi',\psi^{\xi'}_f(N_0))}\cdot X^{p(M,e)}\cdot X^{p(N,f)}\\
      =&\mathrm{LHS}.
    \end{aligned}
   \end{equation*}
  The second equality follows from Lemma \ref{lemp}. The third equality is deduce from the fact $\bigsqcup_c Z^0_{g,c}=\Img\psi^0_g=\bigsqcup_{e+f=g}\Gr_e(FM)\times \Gr_f(FN)$.
  \end{proof}
  \end{proposition}

  Now we need some preparations before we state our first main result.
  \begin{proposition}
  For any $M$, $N$ such that $\Ext^1_{\sC}(M,N)\neq 0$, we have
  \begin{equation*}
    \begin{aligned}
    &\Upsilon(\mathbb{P}\Ext^1_{\sC}(M,N))X_M\cdot X_N\\
    =&\int_{\mathbb{P}\epsilon\in\mathbb{P}\Ext^1_{\sC}(M,N)} g^+(\epsilon,-)\cdot \bar{f}_{M,N}(\epsilon,-)*X_{\rm{mt}\epsilon} +\int_{\mathbb{P}\eta\in\mathbb{P}\Ext^1_{\sC}(N,M)}  \bar{f}_{N,M}(\eta,-)*X_{\rm{mt}\eta}.
    \end{aligned}
    \end{equation*}
  and
  \begin{equation*}
    \begin{aligned}
    &\Upsilon(\mathbb{P}\Ext^1_{\sC}(M,N))X_M\cdot X_N\\
    =&\int_{\mathbb{P}\epsilon\in\mathbb{P}\Ext^1_{\sC}(M,N)} \bar{f}_{M,N}(\epsilon,-)*X_{\rm{mt}\epsilon} +\int_{\mathbb{P}\eta\in\mathbb{P}\Ext^1_{\sC}(N,M)}  g^-(\eta,-)\cdot \bar{f}_{N,M}(\eta,-)*X_{\rm{mt}\eta}.
    \end{aligned}
    \end{equation*}

  \begin{proof}
    We only need to prove the first identity. The proof of the second one is similar. Firstly, we compute the first term on the right hand side.
  \begin{equation*}
    \begin{aligned}
    &\int_{\mathbb{P}\epsilon\in\mathbb{P}\Ext^1_{\sC}(M,N)} g^+(\epsilon,-)\cdot \bar{f}_{M,N}(\epsilon,-)*X_{\rm{mt}\epsilon}\\
    =&\int_{\mathbb{P}\epsilon\in\mathbb{P}\Ext^1_{\sC}(M,N)}\int_g\int_{(M_0,N_0)\in \Img\psi^{\epsilon}_g} t^{g^+(\epsilon,M_0,N_0)}\cdot t^{\Lambda(p(M,[M_0]),p(N,[N_0]))-2c}\cdot t^{2c}\cdot X^{p(\rm mt\epsilon,g)}\\
    =&\int_{\mathbb{P}\epsilon\in\mathbb{P}\Ext^1_{\sC}(M,N)}\int_{e,f}\int_{(M_0,N_0)\in \bigsqcup_{g}\Img\psi_g^{\epsilon}(e,f)} t^{g^+(\epsilon,M_0,N_0)} X^{p(M,e)}\cdot X^{p(N,f)}\\
    =&\int_{e,f}\int_{(M_0,N_0)\in \mathrm{Gr}_e(FM)\times \mathrm{Gr}_f(FN)} \int_{\epsilon\in \Sigma p\Ker\alpha_{M_0,N_0}} t^{2d_{N_0,M_0}}X^{p(M,e)}\cdot X^{p(N,f)}\\
    =&\int_{e,f}\int_{(M_0,N_0)\in \mathrm{Gr}_e(FM)\times \mathrm{Gr}_f(FN)} t^{2d_{N_0,M_0}} \Upsilon(\Sigma p\Ker\alpha_{M_0,N_0})X^{p(M,e)}\cdot X^{p(N,f)}.
    \end{aligned}
  \end{equation*} 
  Next, we compute the second term on the right hand side.
  \begin{equation*}
    \begin{aligned}
    &\int_{\mathbb{P}\eta\in\mathbb{P}\Ext^1_{\sC}(N,M)}  \bar{f}_{N,M}(\eta,-)*X_{\rm{mt}\eta}\\
    =&\int_{\mathbb{P}\eta\in\mathbb{P}\Ext^1_{\sC}(N,M)}\int_g\int_{(N_0,M_0)\in \Img\psi^{\eta}_{g}} t^{\Lambda(p(M,[M_0]),p(N,[N_0]))-2c}\cdot t^{2c}X^{p(\rm mt\eta,g)}\\
    =&\int_{\mathbb{P}\eta\in\mathbb{P}\Ext^1_{\sC}(N,M)}\int_{e,f}\int_{(N_0,M_0)\in \Img\psi_{g}^{\eta}(f,e)} X^{p(M,e)}\cdot X^{p(N,f)}\\
    =&\int_{e,f}\int_{(N_0,M_0)\in \mathrm{Gr}_f(FN)\times \mathrm{Gr}_e(FM)} \int_{\eta\in \P(\Img\alpha'_{N_0,M_0}\cap \Ext^1_{\sC}(N,M))} X^{p(M,e)}\cdot X^{p(N,f)}\\
    =&\int_{e,f}\int_{(N_0,M_0)\in \mathrm{Gr}_f(FN)\times \mathrm{Gr}_e(FM)} \Upsilon(\P(\Img\alpha'_{N_0,M_0}\cap \Ext^1_{\sC}(N,M)))X^{p(M,e)}\cdot X^{p(N,f)}.
    \end{aligned}
  \end{equation*}
  Note that by Lemma \ref{fiberlem}, we have
    $$\dim \Sigma p\Ker\alpha_{M_0,N_0} +\dim (\mathrm{Im}\alpha'_{N_0,M_0}\cap \Ext^1_{\sC}(N,M))=\dim \Ext^1_{\sC}(M,N).$$ 
  Then 
  \begin{equation*}
    \begin{aligned}
    &\Upsilon(\mathbb{P}\Ext^1_{\sC}(M,N))\\
    &=\Upsilon(\P(\Sigma p\Ker\alpha_{M_0,N_0}))+t^{2d'_{M_0,N_0}}\Upsilon(\P(\mathrm{Im}\alpha'_{N_0,M_0}\cap \Ext^1_{\sC}(N,M))),\\
    &=t^{2d_{N_0,M_0}}\Upsilon(\P(\Sigma p\Ker\alpha_{M_0,N_0}))+\Upsilon(\P(\mathrm{Im}\alpha'_{N_0,M_0}\cap \Ext^1_{\sC}(N,M))).
   \end{aligned}
 \end{equation*}
  Therefore, 
   $$\rm{RHS}=\int_{e,f}\int_{(M_0,N_0)\in \mathrm{Gr}_e(FM)\times\mathrm{Gr}_f(FN)} \Upsilon(\mathbb{P}\Ext^1_{\sC}(M,N))X^{p(M,e)}\cdot X^{p(N,f)}=\rm{LHS}.$$
  \end{proof}
  \end{proposition}
  
  The left hand side of the above motivic multiplication formulas do not involve weighted cluster characters. For $\xi\in \Ext^1_{\sC}(M_1,M_2)_M$ and $\xi'\in \Ext^1_{\sC}(N_1,N_2)_N$, let $g$ (resp. $g'$) be the weight function of $(M_1,M_2)$ (resp. $(N_1,N_2)$). We also have a weighted version of multiplication formulas of $g(\xi,-)*X_M$ and $g'(\xi',-)*X_N$, which is a motivic version of \cite[Theorem 4.24]{Chen2021a}. 
  
  \begin{theorem}\label{wtthm}     
    For any $M$, $N$ such that $\Ext^1_{\sC}(M,N)\neq 0$ and any pointwise balanced pair $(h^+,h^-)$, then for any weighted motivic cluster characters $g(\xi,-)*X_M$ and $g'(\xi',-)*X_N$, we have
  \begin{equation}\label{mf1}
    \begin{aligned}
      &\Upsilon(\mathbb{P}\Ext^1_{\sC}(M,N)) (g(\xi,-)*X_M)\cdot (g'(\xi',-)*X_N)\\
        =&\int_{\mathbb{P}\epsilon\in\mathbb{P}\Ext^1_{\sC}(M,N)} h^+(-)\cdot \bar{f}_{M,N}(\epsilon,-)\cdot \T_{g,g'}^{\xi,\xi'}(\epsilon,-)*X_{\rm{mt}\epsilon}+\\
        & \int_{\mathbb{P}\eta\in\mathbb{P}\Ext^1_{\sC}(N,M)}  h^-(-)\cdot \bar{f}_{N,M}(\eta,-)\cdot \hat{\T}^{\xi,\xi'}_{g,g'}(\eta,-)*X_{\rm{mt}\eta}.
    \end{aligned}
  \end{equation}
  \begin{proof}
  Denote by $Z_{e,f}:=\Gr_e(FM)\times \Gr_f(FN)$), $K_{M_0,N_0}:=\P(\Sigma p\Ker\alpha_{M_0,N_0})$ and $K'_{N_0,M_0}:=\P(\Img\alpha'_{N_0,M_0}\cap \Ext^1_{\sC}(N,M))$.By definition, we have that
  \begin{equation*}
    \begin{aligned}
      &\rm{LHS} \\
      =&\Upsilon(\mathbb{P}\Ext^1_{\sC}(M,N))\int_{e,f}\int_{(M_0,N_0)\in Z_{e,f}} t^{g(\xi,\psi_e^{\xi}(M_0))}\cdot t^{g'(\xi',\psi_f^{\xi'}(N_0))}X^{p(M,e)}\cdot X^{p(N,f)}\\
      =&\int_{e,f}\int_{(M_0,N_0)\in Z_{e,f}} t^{g(\xi,\psi_e^{\xi}(M_0))+g'(\xi',\psi_f^{\xi'}(N_0))} \Upsilon(\mathbb{P}\Ext^1_{\sC}(M,N))X^{p(M,e)}\cdot X^{p(N,f)}
    \end{aligned}
  \end{equation*}

  Note that $\T^{\xi,\xi'}_{g,g'}$ and $\hat{\T}_{g,g'}^{\xi,\xi'}$ are independent of choices of $\epsilon\in \Ext^1_{\sC}(M,N)$ and $\eta\in \Ext^1_{\sC}(N,M)$ respectively. Then we have 
  \begin{equation*}
    \begin{aligned}
    &\rm{RHS}\\
    =&\int_{e,f}\int_{(M_0,N_0)\in Z_{e,f}} t^{h^+(M_0,N_0)}\cdot t^{g(\xi,\psi_e^{\xi}(M_0))+g'(\xi',\psi_f^{\xi'}(N_0))}\cdot \Upsilon(K_{M_0,N_0})X^{p(M,e)}\cdot X^{p(N,f)}\\
    +&\int_{e,f}\int_{(M_0,N_0)\in Z_{e,f}} t^{h^-(N_0,M_0)}\cdot t^{g(\xi,\psi_e^{\xi}(M_0))+g'(\xi',\psi_f^{\xi'}(N_0))}\cdot \Upsilon(K'_{N_0,M_0})X^{p(M,e)}\cdot X^{p(N,f)}.
    \end{aligned}
  \end{equation*}
  The property of pointwise balanced pair implies that RHS=LHS. 
  \end{proof}
  \end{theorem}
  
  Notice that the virtual Poincar\'{e} polynomial $\Upsilon(X)$ of $X$ specialized at $t=-1$ is exactly the Euler -Poincar\'{e} characteristic $\chi(X)$. Furthermore, if the skew-symmetric form $\Lambda$ and weight functions $g(\epsilon,-)$ take values in $2\Z$, then the motivic weighted cluster characters $g(\epsilon,-)*X_M$ specialized at $t=-1$ coincides with  the cluster character $CC(M)$ (see Section \ref{sec0}).

  \begin{remark}\label{rem1} (i) If all weight functions involved in Theorem \ref{wtthm} take values in $2\Z$, then the multiplication formula (\ref{mf1}) in the above theorem recovers the cluster multiplication formulas (\ref{in3}) obtained by Palu in \cite{Palu2012}.

  (ii)The multiplication formula (\ref{mf1}) in Theorem \ref{wtthm} is the motivic version of the weighted multiplication formula in \cite[Theorem 4.24]{Chen2021a}.
  \end{remark}

\section{The refined motivic weighted multiplication formula}\label{secref}
  Keller-Plamondon-Qin in \cite{Keller2023} obtained a refined multiplication formula for cluster characters of any 2-Calabi-Yau triangulated categories with a cluster tilting object. let $T$ be a cluster tilting object in $\sC$. Recall the cluster character of $M\in \sC$ is 
     $$CC(M):=CC_T(M)= \sum_e \chi(\rm{Gr}_e(FM)) X^{p(M,e)}.$$
  For any subset $W$ of $\Ext^1_{\sC}(M,N)$, denote by $W_E$ the subset of triangles in $W$ with middle term $E\in \sC$.  For any non-zero vector subspace $V$ of $\Ext^1_{\sC}(M,N)$, they showed a multiplication formula \cite[Theorem 2.10]{Keller2023} as follows
  \begin{equation}\label{KPQ}
    \chi(\P V)CC_T(M)CC_T(N)=\sum_{E\not\cong M\oplus N} (\chi(\P V_E)+ \chi(W_{E}))CC_T(E),
  \end{equation}
  where $W=\P\Ext^1_{\sC}(N,M)\setminus \P(V^{\perp})$, here $V^{\perp}=\Ker(V,-)$ with respect to the non-degenerate bifunctorial form $\beta_{M,N}$ (see Section \ref{sec2.1}).

  The goal of this section is to show that there exists a refined version of the motivic weighted multiplication formula (\ref{mf1}). Fix notations as in Section \ref{motsec}. Remind that the key point of Theorem \ref{wtthm} lies in Lemma \ref{fiberlem}. Let $V$ be a non-zero vector subspace of $\Ext^1_{\sC}(M,N)$. We will reconsider Lemma \ref{fiberlem} from the whole space $\Ext^1_{\sC}(M,N)$ to its subspace $V$. Obviously the fiber of $\psi_g^{\epsilon}$ is independent of choices of $V$. So we have 
  \begin{lemma}\label{reflem}
  \item[(i)] For $\epsilon\in V$, $\mathrm{Gr}_{M_0,N_0}^\epsilon(F(\mathrm{mt}\epsilon))$ is non-empty if and only if 
    $$\epsilon\in \Sigma p\Ker(\alpha_{M_0,N_0})\cap V.$$
  Furthermore in this case, the fiber of $(M_0,N_0)$ satisfies
    $$\mathrm{Gr}_{M_0,N_0}^{\epsilon}(F(\mathrm{mt}\epsilon))\cong \Hom_{B}(M_0,FN/N_0).$$ 
  \item[(ii)] For $\eta\in \mathrm{Ext}^1_{\sC}(N,M)\setminus V^{\perp}$, $\mathrm{Gr}_{N_0,M_0}^\eta(F(\mathrm{mt}\eta))$ is non-empty if and only if 
    $$\eta\in (\mathrm{Im}(\alpha'_{N_0,M_0})\cap \mathrm{Ext}^1_{\sC}(N,M))\setminus (\mathrm{Im}(\alpha'_{N_0,M_0})\cap V^{\perp}).$$
  Furthermore in this case, the fiber of $(N_0,M_0)$ satisfies
     $$\mathrm{Gr}_{N_0,M_0}^{\eta}(F(\mathrm{mt}\eta))\cong \Hom_{B}(N_0,FM/M_0).$$  
  \item[(iii)] Denote by $K_{M_0,N_0}:=\Sigma p\Ker\alpha_{M_0,N_0}$ and $K'_{N_0,M_0}:=\mathrm{Im}\alpha'_{N_0,M_0}\cap \Ext^1_{\sC}(N,M)$. We have the following equation:
   $$\dim (K_{M_0,N_0}\cap V) +\dim K'_{N_0,M_0}/(K'_{N_0,M_0}\cap V^{\perp})=\dim V.$$ 
  \begin{proof}
  The statements $(i)$ and $(ii)$ are clear. We prove the third statement. Since $K_{M_0,N_0}=(K'_{N_0,M_0})^{\perp}$ by Lemma \ref{fiberlem} (iii), then
    \begin{equation*}
    \begin{aligned}
      &\dim (K_{M_0,N_0}\cap V) +\dim (K'_{N_0,M_0}/K'_{N_0,M_0}\cap V^{\perp})\\
      =&\dim (K_{M_0,N_0}\cap V)+\dim (K_{M_0,N_0})^{\perp}\big/ (K_{M_0,N_0}\cup V)^{\perp}\\
      =&\dim (K_{M_0,N_0}\cap V)+\dim(K_{M_0,N_0}\cup V)-\dim(K_{M_0,N_0})=\dim V.
    \end{aligned}
    \end{equation*}
  \end{proof}
  \end{lemma}

  \begin{definition}
  For $0\neq V\subseteq \Ext^1_{\sC}(M,N)$ and a pair $(h^+,h^-)$ of functions on $\Gr(FM)\times \Gr(FN)$, we call $(h^+,h^-)$ is a pointwise balanced pair of $V$, if 
  \begin{align*}
    \Upsilon(\P V)=&t^{h^+(M_0,N_0)}\Upsilon(\P(\Sigma p\Ker\alpha_{M_0,N_0}\cap V))+\\
    &\qquad t^{h^-(N_0,M_0)}\Upsilon(\P(\Img\alpha'_{N_0,M_0}\cap \Ext^1_{\sC}(N,M))\setminus \P(\Img\alpha'_{N_0,M_0}\cap V^{\perp})).
   \end{align*}
  
  \end{definition}

  \begin{example} For $(M_0,N_0)\in \Gr(FM)\times \Gr(FN)$ and $0\neq V\subseteq \Ext^1_{\sC}(M,N)$. Let $g_v(N_0,M_0):=2\dim(\Sigma p\Ker\alpha_{M_0,N_0}\cap V)-2\dim(\Img\alpha'_{N_0,M_0}\cap V^{\perp})$. Then Lemma \ref{reflem}(iii) implies that  $(0,g_v)$ is a pointwise balanced pair of $V$.
  \end{example}

  Thus we can prove the following theorem in the same way as in Theorem \ref{wtthm}. For weight functions $g(\xi,-)$ and $g'(\xi',-)$, recall that $\T^{\xi,\xi'}_{g,g'}$ is the weight function sending $(\epsilon,M_0,N_0)\in MG(M,N)$ to 
      $$g(\xi,\psi_{e}^{\xi}(M_0))+g'(\xi',\psi_{f}^{\xi'}(N_0)),$$ 
  and $\hat{\T}^{\xi,\xi'}_{g,g'}$ is the weight function sending $(\eta,N_0,M_0)\in MG(N,M)$ to 
   $$g(\xi,\psi_{e}^{\xi}(M_0))+g'(\xi',\psi_{f}^{\xi'}(N_0)),$$
  if $(M_0,N_0)\in \rm{Gr}_e(FM)\times \rm{Gr}_f(FN)$.

  \begin{theorem} \label{refthm}    
  For any $M$, $N$ such that $\Ext^1_{\sC}(M,N)\neq 0$ and $V\subseteq \Ext^1_{\sC}(M,N)$ is a non-zero vector subspace, then for any pointwise balanced pair $(h^+,h^-)$ of $V$ and weighted motivic cluster characters $g(\xi,-)*X_M$ and $g'(\xi',-)*X_N$, we have
  \begin{equation}\label{mKPQ}
    \begin{aligned}
      \Upsilon(\mathbb{P}V) g(\xi,-)*X_M &\cdot  g'(\xi',-)*X_N=\int_{\mathbb{P}\epsilon\in\mathbb{P}V}  h^+(-)\cdot \bar{f}_{M,N}(\epsilon,-)\cdot \T_{g,g'}^{\xi,\xi'}(\epsilon,-)*X_{\rm{mt}\epsilon}\\
      &+\int_{\mathbb{P}\eta\in\mathbb{P}\Ext^1_{\sC}(N,M)\setminus \P(V^{\perp})}  h^-(-)\cdot\bar{f}_{N,M}(\eta,-)\cdot \hat{\T}^{\xi,\xi'}_{g,g'}(\eta,-)*X_{\rm{mt}\eta}. 
    \end{aligned}
  \end{equation}
  \begin{proof}
  The proof is similar to Theorem \ref{wtthm}.
  The only difference lies in integral regions. The integral region of the first term in the $\mathrm{RHS}$ becomes 
    $$\int_{\P\epsilon\in \P V}\int_{g,c}\int_{(M_0,N_0)\in Z_{g,c}^{\epsilon}}=\int_{e,f}\int_{(M_0,N_0)\in \Gr_e(FM)\times \Gr_f(FN)}\int_{\P\epsilon\in \P(\Sigma p\Ker\alpha_{M_0,N_0}\cap V)}.$$
  The integral region of the second term in the $\mathrm{RHS}$ becomes 
    \begin{align*}
      &\int_{\P\eta\in \P\Ext^1_{\sC}(N,M)\setminus \P(V^{\perp})}\int_{g,c}\int_{(M_0,N_0)\in Z_{g,c}}\\
      =&\int_{e,f}\int_{(M_0,N_0)\in \Gr_e(FM)\times \Gr_f(FN)}\int_{\P\eta\in \P(K'_{N_0,M_0})\setminus \P(K'_{N_0,M_0}\cap V^{\perp})},
    \end{align*}  
  where $K'_{N_0,M_0}:=\Img\alpha'_{N_0,M_0}\cap \Ext^1_{\sC}(N,M)$. Denote by $K_{M_0,N_0}:=\Sigma p\Ker\alpha_{M_0,N_0}$. Hence, we obtain
  \begin{align*}
    &\rm{RHS}\\
    =&\int_{e,f}\int_{(M_0,N_0)} \big{(}t^{h^+(M_0,N_0)}\Upsilon(\P(K_{M_0,N_0}\cap V))+ t^{h^-(N_0,M_0)}\Upsilon(\P K'_{N_0,M_0}\setminus \P(K'_{N_0,M_0}\cap V^{\perp}))\big{)}\\
    & t^{g(\xi,\psi^{\xi}_e(M_0))}t^{g'(\xi',\psi^{\xi'}_f(N_0))}X^{p(M,e)}\cdot X^{p(N,f)}\\
    =&\rm{LHS}.
  \end{align*}   
  \end{proof}
  \end{theorem}
  
  As a counterpart of Remark \ref{rem1}, we give a remark in refined case.
  \begin{remark}
  Assume that all weight functions and pointwise balanced pair in Theorem \ref{refthm} take values in $2\Z$, then the refined motivic weighted multiplication formula (\ref{mKPQ}) recovers the one (\ref{KPQ}) given by Keller-Plamondon-Qin. 
  \end{remark}

  Although The above refined multiplication formula counts a non-zero vector subspaces of $\Ext^1_{\sC}(M,N)$, it still needs to count elements on the whole space $\Ext^1_{\sC}(N,M)$ of the second term in $\rm{RHS}$. Could we consider a multiplication formula counting subspaces $V$ of $\Ext^1_{\sC}(M,N)$ and $W$ of $\Ext^1_{\sC}(N,M)$? The answer is positive. Since essential point lies in whether the equation
     $$\dim V=\dim(\Sigma p\Ker\alpha_{M_0,N_0}\cap V)+\dim(\Img\alpha'_{N_0,M_0}\cap W)$$
  holds for any $(M_0,N_0)\in \Gr(FM)\times \Gr(FN)$.

  \begin{definition} \label{def3.6}
  (i) Given a non-zero vector subspace $V$ of $\Ext^1_{\sC}(M,N)$, we say a vector subspace $W\subseteq \Ext^1_{\sC}(N,M)$ is a good pair if for any $(M_0,N_0)\in \rm{Gr}_e(FM)\times \rm{Gr}_{f}(FM)$, we have
    $$\dim(V)=\dim(\Sigma p\Ker\alpha_{M_0,N_0}\cap V)+\dim(\Img\alpha'_{N_0,M_0}\cap W).$$
  (ii) Let $(V,W)$ be a good pair, we say $(h^+,h^-)$ is a pointwise balanced pair of $(V,W)$ if 
    $$\Upsilon(\P V)=h^+(M_0,N_0)\Upsilon(\P(\Sigma p\Ker\alpha_{M_0,N_0}\cap V))
    +h^-(N_0,M_0)\Upsilon(\P(\Img\alpha'_{N_0,M_0}\cap W))$$ 
  for any $(M_0,N_0)\in \rm{Gr}_e(FM)\times \rm{Gr}_{f}(FM)$.
  \end{definition}
  
  Thus we can obtain a refined version of motivic weighted multiplication formula with respect to an good pair $(V,W)$  in a similar way in Theorem \ref{refthm}. 
  \begin{theorem}\label{refthm2}     
  For any $M$, $N$ such that $\Ext^1_{\sC}(M,N)\neq 0$,then for any good pair $(V,W)$ with pointwise balanced pair $(h^+,h^-)$ and motivic weighted cluster characters $g(\xi,-)*X_M$ and $g'(\xi',-)*X_N$, we have
  \begin{equation}\label{resformula}
    \begin{aligned}
      &\Upsilon(\mathbb{P}V) (g(\xi,-)*X_M)\cdot (g'(\xi',-)*X_N)\\
      =&\int_{\mathbb{P}\epsilon\in\mathbb{P}V}  h^+(-)\cdot \bar{f}_{M,N}(\epsilon,-)\cdot \T_{g,g'}^{\xi,\xi'}(\epsilon,-)*X_{\rm{mt}\epsilon}\\
      +&\int_{\mathbb{P}\eta\in\mathbb{P}W}  h^-(-)\cdot\bar{f}_{N,M}(\eta,-)\cdot \hat{\T}^{\xi,\xi'}_{g,g'}(\eta,-)*X_{\rm{mt}\eta}.
    \end{aligned}
  \end{equation}
  \end{theorem}
  
  In next section, we will see that there exists such good pair if $\sC$ is a cluster category of an acyclic quiver.

\section{Refined motivic multiplication formulas in hereditary case}\label{SecHer}
  In this section, we want to show that the refined   motivic weighted multiplication formula in Theorem \ref{refthm2} with suitable weight functions $g(\xi,-)$ and $g'(\xi',-)$ and good pair $(V,W)$ recovers the multiplication formula obtained in \cite[Theorem 7.4]{Chen2023}. Now let us fix some notations.

  \subsection{Cluster category of an acyclic quiver}\label{sec5.1}
  For simplicity, let $Q$ be an acyclic quiver with vertex set $\{1,2,\cdots, n\}$. In this subsection, $k$ can be finite field $\F_q$ with $|k|=q$ or $\C$. Let $kQ$ be the path algebra. Denote by $S_i$ the simple $kQ$-module at vertex $i$. The Euler form $\langle-,-\rangle$ on the Grothendieck group $K_0(\mathrm{mod}(kQ))$ is given by  
      $$\langle[M],[N]\rangle=\dim_k \Hom_{kQ}(M,N)-\dim_k \Ext^1_{kQ}(M,N).$$
  Let $\langle-,-\rangle_a$ be the skew-symmetric form of the Euler form. Let $E$ (resp. $B$) be the (resp. skew-symmetric) matrix of $\langle-,-\rangle$ (resp. $\langle-,-\rangle_a$) with respect to $\{[S_i]|\ i=1,\cdots n\}$. Denote by $E^t$ the transpose of $E$, then 
         $$B=E-E^{t}.$$
  We also assume that $B$ is of full rank, which means that there exists a skew-symmetric integral matrix $\Lambda$ such that
       $$-\Lambda\cdot B= I_n,$$
  where $I_n$ is the $n\times n$ identity matrix.
 
  The cluster category $\sC(Q)$ of $Q$ is defined as
     $$\sC(Q):=\sC(kQ):=D^b(\rm{mod}(kQ))/\tau^{-1}[1],$$
  where $[1]=\Sigma$ is the shift functor and $\tau$ is the Auslander-Reiten functor. The cluster category $\sC(Q)$ is a 2-Calabi-Yau category \cite{Buan2006a}. Thus there is a natural non-degenerate bilinear form
     $$(-,-): \Ext^1_{\sC}(X,Y)\otimes \Ext^1_{\sC}(Y,X)\lrw k,$$
  for any $X$, $Y\in \sC(Q)$. Moreover, all isoclasses of indecomposable objects in $\sC(Q)$ can be identified with
       $$\rm{imd}\sC(Q)=\rm{ind}(\rm{mod}(kQ))\cup \{P_1[1],\cdots,P_n[1]\},$$
  where $P_i$ is the indecomposable projective $kQ$-module at vertex $i$. Note that $A':=kQ$ as a $kQ$-module is a cluster tilting object in $\sC(Q)$, there exists an equivalence \cite{Keller2007}
       $$F:=\Hom_{\sC}(A',-): \sC(Q)/A'[1]\lrw \rm{mod}(kQ).$$
  For $\eta\in \Hom_{kQ}(N,\tau M)$, we have the following triangle in $D^b(\rm{mod}(kQ))$
       $$\tau M[-1]\lrw \Ker\eta\oplus \Coker\eta[-1] \lrw N\stackrel{\eta}\lrw \tau M.$$
  Write $\Coker\eta=\tau A\oplus I$ such that $I$ is a injective module and $A$ has no injective direct summand, then we have a triangle in $\sC(Q)$
      $$M\lrw \Ker\eta\oplus A\oplus I[-1] \lrw N\stackrel{\eta}\lrw M[1].$$ 
      
\subsection{Specialized quantum cluster characters}\label{sec4.2}
  In this subsection, $k$ is the finite field $\F_q$.
  Now we can introduce the specialized version of quantum cluster characters given in \cite[Theorem 4.8]{Chen2023}. Denote by $e^*:=E^t\cdot e$, $^*e:=E\cdot e$, and $\underline{\dim} M$ the dimension vector of $M\in \rm{mod}(kQ)$. For $E\in \sC(kQ)$, the (specialized) quantum cluster character of $M\oplus I[-1]$ is
      $$\tilde{X}_{M\oplus I[-1]}:=\int_{e} q^{-\frac{1}{2}\langle e,m-i-e\rangle} |\rm{Gr}_e(M)| X^{-e^*-^*(m-i-e)}.$$
  where $i=\underline{\dim} I$ and $m=\underline{\dim} M$.
   
  Due to Lemma \cite[Lemma 2.3]{Palu2008}, we know that for $\tilde{M}=M\oplus I[-1]\in \sC(kQ)$,
     $$p(\tilde{M},e)=E\cdot i- E\cdot (\dim (FM))+B\cdot e=-e^*-^*(m-i-e).$$
  On the motivic side, for $\tilde{M}\in \sC(\C Q)$, let $\sigma_{\tilde{M}}$ be the following triangle in $\sC(\C Q)$
     $$\tilde{M}\stackrel{id}\lrw \tilde{M}\lrw 0\lrw M[1].$$
  Define a weight function $g_M(-):=g(\sigma_M,-)$ over $MG(0,M)$ as 
     $$g(\sigma_M,-)(\epsilon,M',M'')=
     \begin{cases}
      -\langle e,m-i-e\rangle, &\mbox{if } (M',M'')\in \{0\}\times \Gr_e(F\tilde{M}),\\
      0, &\mbox{eles.}
     \end{cases}$$
  Thus, the motivic weighted cluster character  of $\tilde{M}$  weighted by $g_{\tilde{M}}$ is given by 
  \begin{equation*}
    \begin{aligned} 
    g_{\tilde{M}}*X_{\tilde{M}}&=\int_{e}\int_{(0,M_0)\in {0}\times \rm{Gr}_e(M)} t^{-\langle e,m-i-e\rangle} X^{p(\tilde{M},e)}\\
    &=\int_e\int_{M_0\in \rm{Gr}_e(M)}  t^{-\langle e,m-i-e\rangle}X^{-e^*-^*(m-i-e)}.
    \end{aligned}
  \end{equation*}
  
  \begin{remark}
  Let $M\in \rm{mod}(\Z Q)$ be an indecomposable rigid lattice, set $M^{k}=M\otimes_{\Z} k$. By Qin \cite{Qin2012}, we have $\Upsilon(\Gr_e(M^{\C}))|_{t=q^{\frac{1}{2}}}=|\Gr_e(M^{\F_q})|$, it follows that
      $$(g_M*X_{M^\C})|_{t=q^{\frac{1}{2}}}=\tilde{X}_{M^{\F_q}}.$$
  \end{remark}

\subsection{The refined multiplication formula in hereditary case}
  For any $M$, $N\in \rm{mod}(kQ)$, denote by $[M,N]$ (resp. $[M,N]^1$) the dimension of sapce $\Hom_{kQ}(M,N)$ (resp. $\Ext^1_{kQ}(M,N)$). Let $_D\Hom_{kQ}(N,\tau M)_{\tau A\oplus I}$ be the subset consisting  of $f\in\Hom_{kQ}(N,\tau M)$ such that $\Ker f\cong D$ and $\Coker f\cong \tau A\oplus I$. In \cite[Theorem 7.4]{Chen2023}, Chen-Ding-Zhang have shown the following multiplication formula in the (specialized) quantum torus $\sT_{\Lambda}|_{t=q^{\frac{1}{2}}}$:
  \begin{equation}\label{CDZ}
    \begin{aligned}
     &(q^{[M,N]^1}-1) \tilde{X}_M\cdot \tilde{X}_N= q^{\frac{1}{2}\Lambda(m^*,n^*)}\sum_{E\not\cong M\oplus N} |\Ext^1(M,N)_E|\tilde{X}_E\\
     &+ \sum_{D\not\cong N,A,I}q^{\frac{1}{2}\Lambda((m-a)^*,(n+a)^*)+\frac{1}{2}<m-a,n>}|_D\Hom_{kQ}(N,\tau M)_{\tau A\oplus I}| \tilde{X}_A\cdot\tilde{X}_{D\oplus I[-1]}.
    \end{aligned}
  \end{equation}
  
  Given $\eta\in \Hom_{\C Q}(N,\tau M)$, denote by $D(\eta):=\Ker\eta$, $\Coker\eta:=\tau A(\eta)\oplus I(\eta)$. For simplicity, set $B(\eta):=D(\eta)\oplus I(\eta)[-1]$. The main result in this section is 
  \begin{theorem}\label{thmspe}
  For any $M$, $N\in \rm{mod}(\C Q)$ such that $\Ext^1_{\C Q}(M,N)\neq 0$, we have the following identity.
  \begin{equation}\label{eq2}
    \begin{aligned}
    &\Upsilon(\P\Ext^1_{\C Q}(M,N))(g_M*X_M)\cdot (g_N*X_N)=\int_{\P\epsilon\in \P\Ext^1_{\C Q}(M,N)} t^{\Lambda(m^*,n^*)} g_{\rm{mt}\epsilon}*X_{\rm{mt}\epsilon} \\
    +&\int_{\P\eta \in \P\Hom_{\C Q}(N,\tau M)} t^{\Lambda((m-a(\eta))^*,(n+a(\eta))^*)+\langle m-a(\eta),n\rangle } (g_{A(\eta)}*X_{A(\eta)})\cdot (g_{B(\eta)}*X_{B(\eta)}),
    \end{aligned}
  \end{equation}
  where $a(\eta)=\underline{\dim} \tau^{-1}\Coker\eta$.
  \end{theorem}
  
  To prove the theorem, we need to certify that $(\Ext^1_{\C Q}(M,N),\Hom_{\C Q}(N,\tau M))$ is a good pair in the sense of Definition \ref{def3.6} , then to find a suitable pointwise balanced pair of it. Since
     $$\Ext^1_{\sC(Q)}(M,N)\cong \Ext^1_{\C Q}(M,N)\oplus \Hom_{\C Q}(M,\tau N),$$
  and
     $$\Ext^1_{\sC(Q)}(N,M)\cong \Ext^1_{\C Q}(N,M)\oplus \Hom_{\C Q}(N,\tau M),$$ 
  the non-degenerate bilinear form (mentioned in Section \ref{sec2.1})
  $$\beta_{M,N}: \Ext^1_{\sC(Q)}(M,N)\times \Ext^1_{\sC(Q)}(N,M)\lrw \C$$ 
  can be seen as a direct sum of non-degenerate forms $\beta^{'}_{M,N}$ and $\beta^{''}_{M,N}$ induced by the following isomorphism:
     $$\Ext^1_{\C Q}(M,N)\cong D\Hom_{\C Q}(N,\tau M)\ \ \text{and}\ \  \Hom_{\C Q}(M,\tau N)\cong D\Ext^1_{\C Q}(N,M),$$
  where $D:=\Hom_{\C}(-,\C)$ is the functor taking dual vector spaces. Therefore, the orthogonal complement of $\Ext^1_{\C Q}(M,N)\cap \Sigma p\Ker\alpha_{M_0,N_0}$ in $\Ext^1_{\C Q}(M,N)$ must be $\Hom_{\C Q}(N,\tau M)\cap \Img\alpha'_{N_0,M_0}$. Now we give a rigorous proof of it.

  Let the natural inclusion $i_M:M_0\subseteq M$ and $i_N:N_0\subseteq N$ be the lifts in $\sC(Q)$. Recall (we denote by $\Sigma$ the shift functor and $T:=kQ$ in the following)
  \begin{equation*}
    \begin{aligned}
    \alpha_{M_0,N_0}\colon&\Hom_{\mathcal{C}}(\Sigma^{-1}M,N_0)\oplus\Hom_{\mathcal{C}}(\Sigma^{-1}M,N)\\
    &\longrightarrow\Hom_{\mathcal{C}/(T)}(\Sigma^{-1}M_0,N_0)\oplus\Hom_{\mathcal{C}}(\Sigma^{-1}M_0,N)\oplus\Hom_{\mathcal{C}/(\Sigma T)}(\Sigma^{-1}M,N)\\
    &(a,b) \longmapsto (a\circ\Sigma^{-1}\iota_M,\iota_N\circ a\circ\Sigma^{-1}\iota_M-b\circ\Sigma^{-1}\iota_M,\iota_N\circ a-b)
    \end{aligned}
    \end{equation*}
  and
  \begin{equation*}
    \begin{aligned}
    \alpha'_{N_0,M_0}\colon&\Hom_{\Sigma T}(N_0,\Sigma M_0)\oplus\Hom_{\mathcal{C}}(N,\Sigma M_0)\oplus\Hom_{\Sigma^2T}(N,\Sigma M)\\
    &\longrightarrow\Hom_{\mathcal{C}}(N_0,\Sigma M)\oplus\Hom_{\mathcal{C}}(N,\Sigma M)\\
    &(a,b,c)\longmapsto (\Sigma\iota_M\circ a+c\circ\iota_N+\Sigma\iota_M\circ b\circ\iota_N,-c-\Sigma\iota_M\circ b).
    \end{aligned}
  \end{equation*}

  \begin{lemma}\label{lemdec}
  For $V=\Ext^1_{\C Q}(M,N)\neq 0$, we have $(V,DV)$ is a good pair, where $DV=D\Ext^1_{\C Q}(M,N)\cong \Hom_{\C Q}(N,\tau M)$.
  \begin{proof}
  For any cluster-tilting object $T'$ set $B=\rm{End}_{\sC}(T')$ and $F':=\Hom_{\sC}(T',-)$, we have that
    $$\Hom_{\sC}(X,Y)\cong \Hom_{\sC/(\Sigma T')}(X,Y)\oplus \Hom_{\Sigma T'}(X,Y),$$
  and $\Hom_{\sC/(\Sigma T')}(X,Y)\cong \Hom_B(F'(X),F'(Y))$. Then $\Sigma\alpha_{M_0,N_0}=\alpha_1\oplus \alpha_2$, where $\alpha_1$ is the restriction to $\Hom_{\C Q}(M,\tau N_0)\oplus \Hom_{\C Q}(M,\tau N)$, and $\alpha_2$ is the restriction to $(\Sigma T)(M,\Sigma N_0)\oplus (\Sigma T)(M,\Sigma N)$, which is isomorphic to $\Ext^1_{\C Q}(M,N_0)\oplus \Ext^1_{\C Q}(M,N)$. Similarly, $\alpha'_{N_0,M_0}=\alpha'_1\oplus \alpha'_2$.  Then
    $$\Ker\Sigma\alpha_{M_0,N_0}=\Ker\alpha_1\oplus \Ker\alpha_2\ \ \text{and}\ \ \Img\alpha'_{N_0,M_0}=\Img\alpha'_1\oplus \Img\alpha'_2.$$
  By Lemma \ref{fiberlem}, we have 
    $$(\Sigma p\Ker\alpha_{M_0,N_0})^{\perp}=\Img\alpha'_{N_0,M_0}\cap \Ext^1_{\sC(Q)}(N,M),$$
  it follows that $(p\Ker\alpha_1\oplus p\Ker\alpha_2)^{\perp}=(\Img\alpha'_1\cap \Hom_{\C Q}(N,\tau M))\oplus (\Img\alpha'_2\cap \Ext^1_{\C Q}(N,M))$. Thus the orthogonal complement of $p\Ker\alpha_1$ in $V$ is $\Img\alpha'_1\cap \Hom_{\C Q}(N,\tau M)$=$\Img\alpha'_{N_0,M_0}\cap DV$. Observe that $p\Ker\alpha_2=\Sigma p\Ker\alpha_{M_0,N_0}\cap V$, then we can deduce that 
      $$\dim V=\dim (\Sigma p\Ker\alpha_{M_0,N_0}\cap V)+\dim (\Img\alpha'_{N_0,M_0}\cap DV),$$
    where $DV=D\Ext^1_{\C Q}(M,N)\cong \Hom_{\C Q}(N,\tau M)$.
  \end{proof}     
  \end{lemma}
  
  Next we compute the dimension of $\Sigma p\Ker\alpha_{M_0,N_0}\cap V$.
  \begin{lemma}
   For any $M$, $N\in \rm{mod}(\C Q)$ and any $(M_0,N_0)\in \rm{Gr}_e(M)\times \rm{Gr}_f(N)$, we have 
    $$\dim (\Sigma p\Ker\alpha_{M_0,N_0}\cap V)=[M,N]^1-[M_0,N/N_0]^1,$$
  \begin{proof}
  For $b\in p\Ker\alpha_{M_0,N_0}\cap \Sigma^{-1}V$, we have $F(b)=0$. Applying $F$ to $\alpha_{M_0,N_0}(a,b)$ with $(a,b)\in \Ker\alpha_{M_0,N_0}$, we get
  $F(\iota_N)F(a)=F(b)=0$. Note that $F(\iota_N)=\iota_N$ is injective, it follows that $F(a)=0$, which is equivalent to  $\Sigma a\in \Ext^1_{\C Q}(M,N_0)$. Hence $\Sigma(a\circ\Sigma^{-1}\iota_M)\in \Ext^1_{\C Q}(M_0,N_0)$. Using the duality \cite[Lemma 3.3]{Palu2008}
      $$\Hom_{\Sigma T}(M_0,\Sigma N_0)\cong D\Hom_{\sC/\Sigma T}(N_0,\Sigma M_0)\cong \Ext^1_{\C Q}(M_0,N_0),$$
  then by definition $\Sigma(a\circ\Sigma^{-1}\iota_M)=0$ in $\Hom_{\sC/\Sigma T}(M_0,\Sigma N_0)$. Hence, $(a,b)\in \Ker\alpha_{M_0,N_0}$ with $b\in \Ext^1_{\C Q}(M,N)$ is equivalent to 
    $$\Sigma\iota_N\circ \Sigma a\circ\iota_M= \Sigma b\circ\iota_M.$$
  Notice that $i_M^*:\Ext^1_{\C Q}(M,N_0)\to \Ext^1_{\C Q}(M_0,N_0)$ is surjective, then $b\in \Sigma p\Ker\alpha_{M_0,N_0}\cap \Ext^1_{\C Q}(M,N)$ if and only if $b$ lies in the pullback of $i_M$ and $i_{N_*}$ as follows.
    $$\begin{tikzcd}
     S\arrow[r,"\rm{pr}_1"]\arrow[d,"\rm{pr}_2"] & \Ext^1_{\C Q}(M,N)\arrow[d,"i^*_M"]\\
     \Ext^1_{\C Q}(M_0,N_0)\arrow[r,"i_N*"] &\Ext^1_{\C Q}(M_0,N).
     \end{tikzcd}$$
  Then $\Sigma p\Ker\alpha_{M_0,N_0}\cap \Ext^1_{\C Q}(M,N)=\rm{pr}_1(S)$ and 
  \begin{equation*}
    \begin{aligned}
    \dim \rm{pr}_1(S)=&\dim \Ker i_M^*+\dim \Img i_{N_*}\\
    =& ([M,N]^1-[M_0,N])^1+ ([M_0,N]^1-[M_0,N/N_0]^1)\\
    =& [M,N]^1-[M_0,N/N_0]^1.
    \end{aligned}
  \end{equation*}
  \end{proof}
  \end{lemma}

  \begin{corollary}\label{corpair}
  For any $M$, $N\in \rm{mod}(\C Q)$ such that $V:=\Ext^1_{\C Q}(M,N)\neq 0$. Define 
     $$h^+(M_0,N_0)=2[M_0,N/N_0]^1.$$
  Then $(h^+,0)$ is the pointwise balanced pair of $(V,DV)$ .
  \end{corollary}

  According to Theorem \ref{refthm2}, for good pair $(V,DV)$ with balanced pair $(h^+,0)$ above, the refined weighted multiplication formula of $g_{M}*X_M$ and $g_N*X_N$ is as follows.
  \begin{equation}\label{eq1}
    \begin{aligned}
      &\Upsilon(\mathbb{P}V) (g_M*X_M)\cdot (g_N*X_N)\\
      =&\int_{\mathbb{P}\epsilon\in\mathbb{P}V}  h^+(-)\cdot \bar{f}_{M,N}(\epsilon,-)\cdot \T_{g_M,g_N}^{\sigma_M,\sigma_N}(\epsilon,-)*X_{\rm{mt}\epsilon}\\
      +&\int_{\mathbb{P}\eta\in\mathbb{P}(DV)}  \bar{f}_{N,M}(\eta,-)\cdot \hat{\T}^{\sigma_M,\sigma_N}_{g_M,g_N}(\eta,-)*X_{\rm{mt}\eta}.
    \end{aligned}
  \end{equation}

  Before we prove Theorem \ref{thmspe} using the above equation, We need the following 
  \begin{lemma}[{\cite[Lemma 7.2]{Chen2023}}] 
  we have the following equation:
    $$\Lambda(p(M,e),p(N,f))
    =\Lambda(m^*,n^*)+\langle e,n-f\rangle-\langle f,m-e\rangle.$$   
  \end{lemma}
  As a result, for $(M_0,N_0)\in \rm{Gr}_e(FM)\times \rm{Gr}_f(FN)$, we have 
 \begin{equation*}
    \begin{aligned}
    &h^+(M_0,N_0)+\bar{f}_{M,N}(\epsilon,M_0,N_0)+\T^{\sigma_M,\sigma_N}_{g_M,g_N}(\epsilon,M_0,N_0)+\langle e+f,m+n-e-f\rangle \\
    =&2[M_0,N/N_0]^1-2[M_0,N/N_0]+\Lambda(m^*,n^*)+2\langle e,n-f\rangle\\
    =&\Lambda(m^*,n^*).
  \end{aligned}
  \end{equation*}
  If $\epsilon\in \Ext^1_{kQ}(M,N)_L$, then for $\psi_g^{\epsilon}(L)=(M_0,N_0)\in \rm{Gr}_e(FM)\times \rm{Gr}_f(FN)$, we have $e+f=g$. Put everything together, we reach the corollary below.
  \begin{corollary}\label{first}
  The first term in $\rm{RHS}$ of Equation (\ref{eq1}) coincides with the motivic version of the one of Equation (\ref{eq2}). Namely
  \begin{equation*}
    \begin{aligned}
      &\int_{\epsilon\in \P V} h^+(-)\cdot \bar{f}_{M,N}(\epsilon,-)\cdot \T_{g_M,g_N}^{\sigma_M,\sigma_N}(\epsilon,-)* X_{\rm{mt}\epsilon}\\
      =&\int_{\epsilon\in \P V} h^+(-)\cdot \bar{f}_{M,N}(\epsilon,-)\cdot \T_{g_M,g_N}^{\sigma_M,\sigma_N}(\epsilon,-)\cdot (-g_{\rm{mt}\epsilon})* (g_{\rm{mt}\epsilon}*X_{\rm{mt}\epsilon})\\
      =& \int_{\epsilon\in \P V} t^{\Lambda(m^*,n^*)}g_{\rm{mt}\epsilon}*X_{\rm{mt}\epsilon}.
    \end{aligned}
  \end{equation*}
  \end{corollary}

  Lastly, we compute the second term in $\rm{RHS}$ of Equation (\ref{eq1}).  As mentioned in Section \ref{sec5.1}, for $\eta\in \Hom_{\C Q}(N,\tau M)$, the triangle of $\eta$ in $\sC(Q)$ is as follows
     $$M\stackrel{i}\lrw  A\oplus I[-1]\oplus D \stackrel{p}\lrw N\stackrel{\eta}\lrw M[1],$$
  where $\iota: D:=\Ker\eta\hookrightarrow N$ and $(\tau(p'),r):\tau M\twoheadrightarrow \Coker\eta=:\tau A\oplus I$. Applying $F$ to the above triangle, we have an exact sequence
     $$F(\eta): M\stackrel{(p',0)}\lrw A\oplus D \stackrel{(0,\iota)}\lrw N.$$
  Hence for any $(N_0,M_0)\in \Img\psi^{\eta}_{g}\cap \rm{Gr}_f(N)\times \rm{Gr}_e(M)$, we have $(N_0,p'(M_0))$ belongs to $\rm{Gr}_f(D)\times \rm{Gr}_{e-(m-a)}(A)$.  Conversely, if $(U,V)\in \rm{Gr}_f(D)\times \rm{Gr}_{e-(m-a)}(A)$, then $U\oplus V\subseteq A\oplus D$ and $\psi^{\eta}(U\oplus V)=(V,Fi^{-1}(U))$. Thus, we obtain 
  \begin{lemma}\label{trlem}
  For any $\eta\in \Hom_{\C Q}(N,\tau M)$ and $g$, there is an isomorphism between constructible sets $$\Img\psi^{\eta}_{g}\cap (\rm{Gr}_f(N)\times \rm{Gr}_e(M)) \cong\ \rm{Gr}_f(D)\times \rm{Gr}_{e-(m-a)}(A). $$
  \end{lemma}

  By abuse of notation, we write $\tilde{X}_M:=g_M*X_M$ for $M\in \rm{mod}(\C Q)$. Now we compare $\tilde{X}_{A}\tilde{X}_{D\oplus I[-1]}$ with $X_{D\oplus A\oplus I[-1]}$. Following from Proposition \ref{prop2.11}, we obtain
 
  \begin{lemma}\label{lemcob}
  Denote by $B:=D\oplus I[-1]$. We have the following identity:
    $$\tilde{X}_{A}\cdot\tilde{X}_{B}
    =\bar{f}_{A,B}(0,-)*\T_{g_{A},g_{B}}^{\sigma_{A},\sigma_{B}}(0,-)*X_{A\oplus B},$$
  where $\bar{f}_{A,B}(\epsilon,U,V)=\Lambda(p(A,u),p(B,v))-2[V,A/U]$, and $\T^{\sigma_A,\sigma_B}_{g_A,g_{B}}=-\langle u,a-u \rangle-\langle v,d-i-v \rangle$ for any $(U,V)\in \rm{Gr}_{u}(A)\times \rm{Gr}_v(D)$.
  \end{lemma}

  To compute the power of the second term in $\rm{RHS}$ of Equation $\ref{eq2}$, we also need another lemma.

  \begin{lemma}[{\cite[Lemma 7.3]{Chen2023}}] For $\eta\in \Hom_{\C Q}(N,\tau M)$ and fix notations as above. Denote $i=\underline{\dim} I(\eta)$, $d=\underline{\dim} D(\eta)$ and $a=\underline{\dim} A(\eta)$, the we have
    $$\Lambda((m-a)^*,(n+a)^*)+\langle m-a,n\rangle +\Lambda(a^*,(d-i)^*)=\Lambda(m^*,n^*)+\langle m-a,n-a\rangle.$$
  \end{lemma}
  \ \ \\
  Now we can compute $t^{\Lambda((m-a)^*,(n+a)^*)+\langle m-a,n\rangle}\tilde{X}_A\tilde{X}_{D\oplus I[-1]}$ for $0\neq \eta\in \Hom_{\C Q}(N,\tau M)$. For $(U,V)\in \rm{Gr}_{u}(A)\times \rm{Gr}_v(D)$, set $\bar{f}(U,V):=$ $\bar{f}_{A,D\oplus I[-1]}(0,(U,V))$ and $\T(U,V):=\T^{\sigma_A,\sigma_{D\oplus I[-1]}}_{g_A,g_{D\oplus I[-1]}}(0,U,V)$, then we have that
  \begin{equation*}
  \begin{aligned}
    &\Lambda((m-a)^*,(n+a)^*)+<m-a,n>+  \bar{f}(U,V)+\T(U,V)\\
    =&\Lambda((m-a)^*,(n+a)^*)+<m-a,n>+\Lambda(p(A,U),p(D\oplus I[-1],V))\\
    &\ -\langle u,a-u \rangle-\langle v,d-i-v \rangle -2[V,A/U]\\
    =&\Lambda((m-a)^*,(n+a)^*)+<m-a,n>+\Lambda(a^*,(d-i)^*)+2\langle u,d-i-v \rangle\\
    &\  -<u+v,a+d-i-u-v>-2[V,A/U] \\
    =&\Lambda(m^*,n^*)+\langle m-a,n-a\rangle+2<u,d-i-v>-<u+v,a+d-i-u-v>\\
    &\ -2[V,A/U].\\
  \end{aligned}
  \end{equation*}

 Using Lemma \ref{trlem}, set $M_0=Fi^{-1}(U)$, and $N_0:=V$, we have $f=v$, $u=e-(m-a)$ and $M/M_0\cong A/U$. Then  
  \begin{equation}
  \begin{aligned}
    &\Lambda((m-a)^*,(n+a)^*)+<m-a,n>+  \bar{f}(U,V)+\T(U,V)\\
    =& \Lambda(m^*,n^*)+\langle m-a,n-a\rangle+2<e-(m-a),d-i-f>\\
    &\ -<e+f-(m-a),a+d-i-(e+f-m+a)>-2[N_0,M/M_0]\\
    =&\Lambda(m^*,n^*)+2\langle e,n-f \rangle-\langle e+f,m+n-e-f \rangle-2[N_0,M/M_0]\\
    =& \bar{f}_{N,M}(\eta,N_0,M_0)\cdot \hat{\T}_{g_M,g_N}(\eta,N_0,M_0).
  \end{aligned}
  \end{equation}
  The second equality comes from $d-i=n-\tau m +\tau a$,  and the fact that $\langle [N'],\tau m-\tau a \rangle=-\langle m-a, [N'] \rangle$ for any $N'\in \rm{mod}(\C Q)$ if $A\not\cong M$ (because $\rm{mod}(\C Q)$ is hereditary).
  
  Combining last equation with Lemma \ref{lemcob}, we can deduce the following corollary.
  \begin{corollary}\label{second}
  The second term in $\rm{RHS}$ of Equation (\ref{eq1}) coincides with the one of Equation (\ref{eq2}). Namely,
  \begin{equation*}
  \begin{aligned}
    &\int_{\mathbb{P}\eta\in\mathbb{P}\Hom_{kQ}(N,\tau M)} \bar{f}_{N,M}(\eta,-)\cdot \hat{\T}^{\sigma_M,\sigma_N}_{g_M,g_N}(\eta,-)*X_{\rm{mt}\eta}\\
    =&\int_{\P\eta\in \P\Hom_{kQ}(N,\tau M)} t^{\Lambda((m-a(\eta))^*,(n+a(\eta))^*)+<m-a(\eta),n>}*\\
    &\qquad \qquad \bar{f}_{A(\eta),D(\eta)\oplus I(\eta)[-1]}(0,-)*\T_{A(\eta),D(\eta)\oplus I(\eta)[-1]}(0,-)* X_{A(\eta)\oplus D(\eta)\oplus I(\eta)[-1]}\\
    =&\int_{\P\eta\in \P\Hom_{kQ}(N,\tau M)} t^{\Lambda((m-a(\eta))^*,(n+a(\eta))^*)+\frac{1}{2}<m-a(\eta),n>} \tilde{X}_{A(\eta)}\cdot\tilde{X}_{D(\eta)\oplus I(\eta)[-1]}.
  \end{aligned}
  \end{equation*}
  where $\T_{A(\eta),D(\eta)\oplus I(\eta)[-1]}:=\T^{\sigma_A,\sigma_{D\oplus I[-1]}}_{g_{A(\eta)},g_{D(\eta)\oplus I(\eta)[-1]}}$.
  \end{corollary}
  In the end, we have proved Theorem \ref{thmspe} by Corollary \ref{first} and Corollary \ref{second}.

\subsection{A refined multiplication formula for V of dimension 1}
  In the previous subsection, we have shown a refined motivic multiplication formula for $0\neq V\subseteq \Ext^1_{\C Q}(M,N)$. In this subsection, we want to show another refined motivic multiplication formula for $V\subseteq \Ext^1_{\C Q}(M,N)$ of dimension 1. Recall that the non-degenerate form $\beta_{M,N}$ defining the 2-Calabi-Yau property on $\sC(Q)$ is the direct sum of non-degenerate forms $\beta'_{M,N}:\Ext^1_{\C Q}(M,N)\times \Hom_{\C Q}(N,\tau M)\to \C$ and $\beta^{''}_{M,N}:\Hom_{\C Q}(M,\tau N)\times \Ext^1_{\C Q}(N,M)\to \C$. Let $V$ be an subspace of $\Ext^1_{\C Q}(M,N)$ of dimension 1. Then 
     $$V^{\perp}=V'\oplus \Ext^1_{\C Q}(N,M),$$ 
  where $V'$ is the orthogonal complement of V in $\Hom_{\C Q}(N,\tau M)$ with respect to $\beta'_{M,N}$. From the proof of Lemma \ref{lemdec}, it follows that $\Img\alpha'_{N_0,M_0}=\Img\alpha'_1\oplus \Img \alpha'_2$. Then for any $(M_0,N_0)\in \Gr_e(FM)\times \Gr_f(FN)$, by the definition of $\alpha'_i$, we have
    $$\Ext^1_{\sC}(N,M)\cap \Img \alpha'_{N_0,M_0}=\Hom_{\C Q}(N,\tau M)\cap (\Img\alpha'_{N_0,M_0})\oplus (\Ext^1_{\C Q}(N,M)\cap \Img\alpha'_{N_0,M_0}). $$
  and
     $$V^{\perp}\cap \Img \alpha'_{N_0,M_0}=V'\cap (\Img\alpha'_{N_0,M_0})\oplus (\Ext^1_{\C Q}(N,M)\cap \Img\alpha'_{N_0,M_0}). $$
  Therefore we obtain 
  \begin{equation*}
  \begin{aligned}
    U_{N_0,M_0}:=&(\Ext^1_{\sC(Q)}(N,M) \cap \Img\alpha'_{N_0,M_0})/ (V^{\perp}\cap \Img \alpha'_{N_0,M_0})\\
    \cong & (\Hom_{\C Q}(N,\tau M)\cap \Img\alpha'_{N_0,M_0})/(V'\cap \Img\alpha'_{N_0,M_0})=: U'_{N_0,M_0}.
  \end{aligned}
  \end{equation*}
  As a consequence, by the Lemma \ref{reflem} (iii), we know that
    $$\dim(\Sigma p\Ker\alpha_{M_0,N_0}\cap V)+\dim (U'_{N_0,M_0})=1.$$
  Notice that $F_{N_0,M_0}(V):=\P(\Hom_{\C Q}(N,\tau M) \cap \Img\alpha'_{N_0,M_0})\setminus \P(V'\cap \Img \alpha'_{N_0,M_0})$ is an affine fibration over $\P U'_{N_0,M_0}$ of rank $d_V(N_0,M_0):=\dim (V^{\perp}\cap \Img \alpha'_{N_0,M_0})$, it follows that
     $$\Upsilon(\Sigma p\Ker\alpha_{M_0,N_0}\cap V)+t^{-2d_V(N_0,M_0)}\Upsilon(F_{N_0,M_0}(V))=1.$$
  Finally we deduce that $(0,-2d_V)$ is a pointwise balanced pair of $V$.

  Denote by $m:=\underline{\dim} M$ and $n:=\underline{\dim} N$. Define a weight function $f_V$ by setting 
    $$f_V(\epsilon,M_0,N_0)=-2[M_0,N/N_0]+\Lambda(m^*,n^*)+2\langle e,n-f\rangle$$ 
  for any $\epsilon\in V$ and $(M_0,N_0)\in \Gr_e(FM)\times \Gr_f(FN)$.

  \begin{proposition}
    For any $M$, $N\in \rm{mod}(\C Q)$ such that $\Ext^1_{\C Q}(M,N)\neq 0$ and for $0\neq \epsilon\in \Ext^1_{\C Q}(M,N)$, we have the following identity.
    \begin{equation}
      \begin{aligned}
      &(g_M*X_M)\cdot (g_N*X_N)=  f_V(\epsilon,-)\cdot (g_{\rm{mt}\epsilon}*X_{\rm{mt}\epsilon}) \\
      +&\int_{\P\eta \in \P\Hom_{\C Q}(N,\tau M)\setminus \P V'}  (-2d_V)\cdot \bar{f}_{N,M}(\eta,-)\cdot \hat{\T}^{\sigma_M,\sigma_N}_{g_M,g_N}(\eta,-)*X_{\rm{mt}\eta},
      \end{aligned}
    \end{equation}
  where $V$ is the subspace of $\Ext^1_{\C Q}(M,N)$ spanned by $\epsilon$ and $V'\subseteq \Hom_{\C Q}(N,\tau M)$ is the orthogonal complement of $V$ with respect to $\beta'_{M,N}$. 
  \begin{proof}
  Applying Theorem \ref{refthm} to the above subspace $V$ and the pointwise balanced pair $(0,-2d_V)$, then the refined motivic multiplication formula for $g_M*X_M$ and $g_N*X_N$ is precisely as above.
  \end{proof}
  \end{proposition}

\subsection{Cluster multiplication formulas with initial cluster characters}
  Note that $\bigoplus_{1\leq i\leq n}P_i[1]$ is a cluster-tilting object in $\sC(\C Q)$, which is always chosen to be the initial cluster-tilting object. In Section \ref{sec4.2}, we have defined a motivic weighted cluster character of $\tilde{M}=M\oplus I[-1]$ for any $M\in \rm{mod}Q$ and injective module $I$. Notice that any $\tilde{M}=M\oplus I[-1]\cong M\oplus P[1]$, where $P=\nu^{-1}(I)$ projective, and $\nu=D\Hom_{\C Q}(-,A)$ is the Nakayama functor. In particular, the AR translation functor $\tau=\nu[-1]$. Now we give another motivic weighted cluster character of $\tilde{M}$ by the weight function $\tilde{g}_{\tilde{M}}(\sigma_{\tilde{M}},-)$, defined as
    $$\tilde{g}_{\tilde{M}}(\sigma_{\tilde{M}},M',M'')=
    \begin{cases}
    \langle p-e,m-e \rangle,  &\mbox{if } (M',M'')\in \{0\}\times \Gr_e(F\tilde{M}),\\
    0, &\mbox{eles.}
    \end{cases}$$
  Thus, the motivic weighted cluster character  of $\tilde{M}$  weighted by $g_{\tilde{M}}$ is given by 
  \begin{equation*}
    \begin{aligned} 
    \tilde{g}_{\tilde{M}}*X_{\tilde{M}}&=\int_{e}\int_{(0,M_0)\in {0}\times \rm{Gr}_e(M)} t^{\langle p-e,m-e\rangle} X^{p(\tilde{M},e)}\\
    &=\int_e\int_{M_0\in \rm{Gr}_e(M)}  t^{\langle p-e,m-e\rangle}X^{(p-e)^*-^*(m-e)}.
    \end{aligned}
  \end{equation*}
  Note that $p^*=(\dim P)^*=^*(\dim \nu(P))=^*i$, then $p(\tilde{M},e)=-e^*-^*(m-i-e)=(p-e)^*-^*(m-e)$.

  Our goal in this section is to prove the following theorem.

  \begin{theorem}\label{Thm4.13}
  For any $M\in \rm{mod}(\C Q)$ and injective $\C Q$-module $I$, we have
  \begin{equation}\label{eq4.6}
    \begin{aligned}
    &(t^{[M,I]}-1)(g_{M}*X_M)\cdot X_{I[-1]} \\
    = &t^{\Lambda(^*i,^*m)}\big(\sum_{[M'],[I']\not\cong I} \Upsilon(_{M'}\Hom_{\C Q}(M, I)_{I'}) g_{M'\oplus I'[-1]}*X_{M'\oplus I'[-1]}\\
      & \ \ +t^{\langle m,i \rangle } \sum_{[M''],[P']\not\cong P} \Upsilon(_{P'}\Hom_{\C Q}(P, M)_{M''}) \tilde{g}_{M''\oplus P'[1]}*X_{M''\oplus P'[1]}\big),
    \end{aligned}
  \end{equation}
  and 
  \begin{equation}\label{eq4.7}
    \begin{aligned}
    &(t^{[M,I]}-1) X_{I[-1]}\cdot (g_{M}*X_M) \\
    = &t^{\Lambda(^*m,^*i)}\big(\sum_{[M'],[I']\not\cong I} \Upsilon(_{M'}\Hom_{\C Q}(M, I)_{I'}) g_{M'\oplus I'[-1]}*X_{M'\oplus I'[-1]}\\
      & \ \ +t^{-\langle m,i \rangle } \sum_{[M''],[P']\not\cong P} \Upsilon(_{P'}\Hom_{\C Q}(P, M)_{M''}) \tilde{g}_{M''\oplus P'[1]}*X_{M''\oplus P'[1]}\big),
    \end{aligned}
  \end{equation}
  where $P=\nu^{-1}(I)$.
  \end{theorem}

  Although the left hand side of Equation (\ref{eq4.6}) counts $\Upsilon(\P\Hom_{\C Q}(M,I))$ instead of $\Ext^1_{\sC(\C Q)}(M,I[-1])$, we do not need to use any refined version. Indeed, we have
    $$\Ext^1_{\sC(\C Q)}(M,I[-1])\cong \Hom_{\C Q}(M,I).$$
  Applying Theorem \ref{wtthm} to $g_{M}*X_M$, $X_{I[-1]}$ and some balanced pair $(h^+,h^-)$, then we have
  \begin{equation*}
    \begin{aligned}
    &\Upsilon(\P \Hom_{\C Q}(M,I))(g_{M}*X_M)\cdot X_{I[-1]}\\
    =&\int_{\epsilon\in \P \Hom_{\C Q}(M,I) } h^+(-)\cdot \bar{f}_{M,I[-1]}(\epsilon,-)\cdot \T_{g_M,0}^{\sigma_M,\sigma_{I[-1]}}(\epsilon,-)*X_{\rm{mt}\epsilon}\\
    +&\int_{\eta\in \P \Hom_{\C Q}(P[1],M[1]) } h^-(-)\cdot \bar{f}_{P[1],M}(\eta,-)\cdot \hat{\T}_{0,g_M}^{\sigma_{P[1]},\sigma_{M}}(\eta,-)*X_{\rm{mt}\eta}.
    \end{aligned}
  \end{equation*}

  Firstly, we compute the exponents of the two terms in the right hand side.
  
  \begin{lemma}\label{lem1}
  For any $M_0\in \rm{Gr}_e(M)$ and $0\neq \epsilon\in \Ext^1_{\sC}(M,I[-1])$ such that $(M_0,0)\in \Img\psi^{\epsilon}$. Write $\rm{mt}\epsilon=M'\oplus I'[-1]$, then we have
  \begin{equation*}
    \bar{f}_{M,I[-1]}(\epsilon,M_0,0)+\T_{g_M,0}^{\sigma_M,\sigma_{I[-1]}}(\epsilon,M_0,0) -g_{\rm{mt}\epsilon}(\sigma_{\rm{mt}\epsilon},0,L_0)=\Lambda(^*i,^*m),
  \end{equation*}
  where $L_0\in (\psi^{\epsilon})^{-1}(M_0,0)$.
  \begin{proof} Note that the fiber of $(\psi^{\epsilon})^{-1}(M_0,0)$ is $M_0$ itself, then 
    \begin{equation*}
      \begin{aligned}
      &\bar{f}_{M,I[-1]}(\epsilon,M_0,0)+\T_{g_M,0}^{\sigma_M,\sigma_{I[-1]}}(\epsilon,M_0,0) -g_{\rm{mt}\epsilon}(\sigma_{\rm{mt}\epsilon},0,M_0)\\
      =& \Lambda(-e^*-^*(m-e),^*i)-\langle e,m-e \rangle +\langle e, m'-i'-e \rangle \\
      =& -\Lambda (^*m,^*i)+ \Lambda(B\cdot e,^*i)-\langle e,m-e \rangle +\langle e,m-i-e \rangle\\
      =&\Lambda (^*i,^*m).
      \end{aligned}
     \end{equation*}
     where $i'=\dim I'$ and $m'=\dim M'$. In particular $m-i=m'-i'$.
  \end{proof}
  
  \end{lemma}

  \begin{lemma}\label{lem2}
  For any $M_0\in \rm{Gr}_e(M)$ and $0\neq \eta\in \Ext^1_{\sC}(P[1],M)$ such that $(0,M_0)\in \Img\psi^{\eta}$. Write $\rm{mt}\eta=M''\oplus P'[1]$, then we have
  \begin{equation*}
    \bar{f}_{P[1],M}(\epsilon,0,M_0)+\hat{\T}_{0,g_M}^{\sigma_{P[1]},\sigma_M}(\epsilon,0,M_0) -\tilde{g}_{\rm{mt}\eta}(\sigma_{\rm{mt}\eta},0,L'_0)=\Lambda (^*i,^*m)-\langle m-2e,i \rangle,
  \end{equation*}
  where $L'_0\in (\psi^{\eta})^{-1}(0,M_0)$.
  \begin{proof}
  Since fiber space $[0,M/M_0]$ at $(0,M_0)$ is just one point space and $\dim L'_0=e-\dim \Img\eta[-1]=:e'$, then  
  \begin{equation*}
    \begin{aligned}
      &\bar{f}_{P[1],M}(\epsilon,0,M_0)+\hat{\T}_{0,g_M}^{\sigma_{P[1]},\sigma_M}(\epsilon,0,M_0) -\tilde{g}_{\rm{mt}\eta}(\sigma_{\rm{mt}\eta},0,L'_0)\\
      =& \Lambda(-e^*-^*(m-e),^*i)-\langle e,m-e \rangle -\langle p'-e', m''-e' \rangle \\
      =& -\Lambda (^*m,^*i)+ \Lambda(B\cdot e,^*i)-\langle e,m-e \rangle -\langle p-e,m-e \rangle\\
      =&\Lambda (^*i,^*m)+\langle e,i \rangle-\langle p,m-e\rangle \\
      =& \Lambda (^*i,^*m)-\langle m-2e,i \rangle.
    \end{aligned}
  \end{equation*}
  where $p'=\dim P'$ and $m''=\dim M''$. The last equality follows from $\langle p,m-e \rangle=\langle m-e,i \rangle$.
  \end{proof}
  \end{lemma}

  Next we want to find a suitable balanced pair. Let us compute the dimension of $\Img\alpha'_{0,M_0}\cap \Ext^1_{\sC}(P[1],M)$, which is a good candidate for balanced pairs.

  \begin{lemma}
  For any $M_0\in \Gr_e(M)$, we have 
    $$\dim (\Img\alpha'_{0,M_0}\cap \Ext^1_{\sC}(P[1],M))= [P,M_0]=\langle p,e \rangle.$$
   \begin{proof}
  By Lemma \ref{fiberlem}, for any $M_0\in \Gr_e(M)$, $\eta\in \Img\alpha'_{0,M_0}\cap \Ext^1_{\sC}(P[1],M)$ if and only if $(0,M_0)\in \Img\psi^{\eta}$. Consider $\eta\in \Hom_{\C Q}(P[1],M[1])$ with middle term $M''\oplus P'[1]$, the $F\eta$ is as follows
      $$F\eta: \Img\eta[-1]\hookrightarrow M\stackrel{Fi}\lrw M''\lrw 0.$$
  Note that $\Ker(Fi)=\Img\eta[-1]$, then $M_0= Fi^{-1}(L_0)$ for some $L_0\subseteq M''$ if and only if $\Img\eta[-1]\subseteq M_0$. The last condition is equivalent to that $\eta[-1]$ lies in the image of 
     $$(i_{M_0})_*: \Hom_{\C Q}(P,M_0)\lrw \Hom_{\C Q}(P,M).$$
  In particular $(i_{M_0})_*$ is injective. Hence 
     $$\Img\alpha'_{0,M_0}\cap \Ext^1_{\sC}(P[1],M)= \Img (i_{M_0})_*[1],$$
  and $\dim (\Img\alpha'_{0,M_0}\cap \Ext^1_{\sC}(P[1],M))=[P,M_0]=\langle p,e \rangle$.
  
  \end{proof}
  \end{lemma}

  \begin{corollary}\label{cor3}
  For any $M_0\in \Gr_e(M)$, we have 
    $$\dim p\Sigma\Ker\alpha_{M_0,0}=\langle p,m-e \rangle=\langle m-e,i \rangle.$$
  \end{corollary}

  \subsection*{Proof of Theorem 4.13}
  \begin{proof}
  Take a balanced pair as $(0, 2\langle m-e,i \rangle)$ for each $(M_0,0)\in \Gr_e(M)\times \{0\}$. By Lemma \ref{lem1} and \ref{lem2} and Corollary \ref{cor3}, we obtain 
  \begin{equation*}
    \begin{aligned}
    &\Upsilon(\P \Hom_{\C Q}(M,I))(g_{M}*X_M)\cdot X_{I[-1]}\\
    =&\int_{\epsilon\in \P \Hom_{\C Q}(M,I) } t^{\Lambda(^*i,^*m)}g_{\rm{mt}\epsilon}*X_{\rm{mt}\epsilon}
    +\int_{\eta\in \P \Hom_{\C Q}(P[1],M[1]) } t^{\Lambda(^*i,^*m)+\langle m,i \rangle}\tilde{g}_{\rm{mt}\eta}*X_{\rm{mt}\eta}\\
    =&\text{LHS of Equation }\ref{eq4.6}.
    \end{aligned}
  \end{equation*}
  The proof of Equation (\ref{eq4.7}) is  similar by taking the balanced pair as $(0,2\langle e,i \rangle)$ for each $(0,M_0)\in \{0\} \times \Gr_e(M)$.
  \end{proof}

\subsection*{Acknowledgement}{All authors are partially supported by Natural Science Foundation of China
[Grant No. 12031007].}

\bibliographystyle{plain}
\bibliography{myref}

\end{document}